\documentclass[lettersize,journal]{IEEEtran}
\usepackage{amsmath,amsfonts}
\usepackage{algorithmic}
\usepackage{algorithm}
\usepackage{array}
\usepackage[caption=false,font=normalsize,labelfont=sf,textfont=sf]{subfig}
\usepackage{textcomp}
\usepackage{stfloats}
\usepackage{url}
\usepackage{verbatim}
\usepackage{graphicx}
\usepackage{cite}
\usepackage{mathtools}
\usepackage{tikz}
\usepackage{multirow}
\usepackage{threeparttable}
\usepackage{booktabs}       

\usepackage{amsthm}
\usepackage{amssymb}
\usepackage{mathtools}
\usepackage{hyperref}
\usepackage{color}
\usepackage{float}
\usepackage{placeins}
\usepackage{csquotes}
\usepackage{makecell}
\usepackage{tabularx}
\usepackage{lineno}


\theoremstyle{definition}
\newtheorem{theorem}{Theorem} 

\theoremstyle{definition}
\newtheorem{lemma}{Lemma} 

\theoremstyle{definition}
\newtheorem{proposition}{Proposition}

\theoremstyle{definition}

\theoremstyle{remark}

\theoremstyle{assumption}
\newtheorem{assumption}{Assumption}

\begin{document}

\title{A Learning-Based Inexact ADMM for Solving Quadratic Programs}

\author{Xi Gao, Jinxin Xiong, Linxin Yang, Akang Wang, Weiwei Xu, and Jiang Xue~\IEEEmembership{Senior Member,~IEEE}%
\thanks{\textit{Corresponding authors: Jiang Xue; Akang Wang.}}%
\thanks{Xi Gao and Jiang Xue are with the School of Mathematics and Statistics, Xi'an Jiaotong University, Xi'an, China (e-mail: \mbox{gaoxi@stu.xjtu.edu.cn}; \mbox{x.jiang@xjtu.edu.cn}).}%
\thanks{Jinxin Xiong, Linxin Yang, and Akang Wang are with Shenzhen International Center for Industrial and Applied Mathematics, Shenzhen Research Institute of Big Data, Shenzhen Loop Area Institute, and also with the School of Data Science, The Chinese University of Hong Kong, Shenzhen, China (e-mail: \mbox{jinxinxiong@link.cuhk.edu.cn}; \mbox{linxinyang@link.cuhk.edu.cn}; \mbox{wangakang@sribd.cn}).}%
\thanks{Weiwei Xu is with the School of Mathematics and Statistics, Nanjing University of Information Science and Technology, Nanjing, China (e-mail: \mbox{wwxu@nuist.edu.cn}).}%
}

\maketitle

\begin{abstract}
Convex quadratic programs (QPs) serve as a fundamental computational primitive across diverse domains, including financial optimization, control systems, and machine learning. 
The alternating direction method of multipliers (ADMM) has emerged as a leading first-order approach due to its iteration efficiency, as exemplified by the state-of-the-art OSQP solver. 
Recently, machine learning-enhanced optimization algorithms have shown significant promise in accelerating the solution process.
In this work, we introduce a neural-accelerated ADMM variant that replaces exact subproblem solutions with learned approximations generated by a Long Short-Term Memory network. 
{A safeguard mechanism is incorporated to revert to conventional optimization steps when the approximate solutions fail to satisfy prescribed inexactness conditions. Based on this, we establish the convergence of the proposed learning-based inexact ADMM framework.
Extensive experimental results on convex QPs demonstrate that our method achieves superior solution accuracy compared to existing learning-based approaches, while delivering speedups of up to $10\times$, $17\times$, and $14\times$ over the traditional solvers Gurobi, SCS and OSQP, respectively, with an objective gap as low as 0.05\%. 
When directly extended to nonconvex problems with linear constraints, our method achieves speedups of up to $200\times$ over IPOPT, while empirically attaining high-quality stationary solutions with an objective gap of about 0.5\% relative to IPOPT.
}

\end{abstract}


\section{Introduction}\label{sec:Intro}
\IEEEPARstart{Q}{uadratic} \textit{Programs} (QPs), as a fundamental class of constrained optimization problems, find widespread applications across diverse domains, ranging from finance and engineering to machine learning. In finance, a prominent example is portfolio optimization, where QPs are employed to allocate assets by optimizing the trade-off between expected returns and risk exposure~\cite{boyd2013performance, xiu2023fast, shi2022cardinality}. 
Control engineering, on the other hand, frequently relies on solving sequential QPs at each time step, necessitating extremely rapid computation to satisfy real-time control requirements~\cite{garcia1989model, borrelli2017predictive, adeoye2024inexact, liufu2025adaptive}. Similarly, in machine learning, numerous unconstrained optimization problems—such as LASSO regression~\cite{tibshirani1996regression}, support vector machines~\cite{cortes1995support}, and Huber regression~\cite{huber2011robust}—can be equivalently reformulated as QPs, as demonstrated in~\cite{stellato2020osqp}. 
Given this extensive applicability, the development of efficient and accurate QP-solving algorithms remains a critical research pursuit.

\textit{Interior-point methods} (IPMs) constitute a fundamental class of algorithms for solving QPs and serve as the computational backbone of leading commercial solvers including MOSEK~\cite{aps2019mosek} and Gurobi~\cite{gurobi2023}. 
As a second-order method, IPMs have gained widespread adoption for various nonlinear programs~\cite{wachter2006implementation}. 
In practice, modern implementation predominantly employs an infeasible path-following primal-dual approach, which initiates from potentially infeasible points and iteratively solves the perturbed KKT conditions via Newton's method, while enforcing centrality conditions via line search. 
These methods generate iterates that follow central paths, progressively reducing infeasibility while converging to optimal solutions. 
Each iteration requires solving a Newton system through matrix factorization, with computational complexity up to $\mathcal{O}(n^3)$ for dense matrices. 
Consequently, despite exhibiting quadratic convergence rates, IPMs face significant challenges in scaling to large-scale or ill-conditioned problems. Furthermore, the central path trajectory inherently limits the effectiveness of warm-start strategies~\cite{john2008implementation}, presenting additional computational bottlenecks for practical implementation.

\textit{First-order methods}, which rely solely on gradient information, have emerged as powerful alternatives to second-order methods in large-scale optimization due to their computational efficiency. 
Notable examples include the \textit{primal-dual hybrid gradient} (PDHG~\cite{he2014convergence}) algorithm and the \textit{alternating direction method of multipliers} (ADMM~\cite{boyd2011distributed}).
The work of~\cite{applegate2021practical} pioneered the application of PDHG to linear programs, achieving accuracy comparable to traditional solvers while demonstrating enhanced computational efficiency through matrix-free operations. 
Efficient parallelization is enabled by the algorithmic structure, allowing significant time reduction through GPU acceleration.
Building upon these advancements, the work of~\cite{lu2023practical} extended the PDHG framework to QPs, attaining linear convergence under mild regularity conditions while maintaining hardware-agnostic implementation across both GPU and CPU architectures.
Mathematically equivalent to \textit{Douglas-Rachford splitting}~\cite{lions1979splitting}, ADMM is particularly suitable for noise-prone environments where computationally efficient, near-optimal solutions are acceptable.
The \textit{splitting conic solver} (SCS)~\cite{o2016conic}, a leading ADMM-based implementation, solves convex conic programs through careful dual variable tuning. 
While any QP can be reformulated as a conic program, such reformulation is often inefficient from a computational point of view.
Our focus, the \textit{operator splitting quadratic programming} (OSQP~\cite{stellato2020osqp}) solver, represents a specialized ADMM variant optimized for convex QPs through adaptive parameter updates that alleviate the computational burden of expensive matrix factorization.

Recently, \textit{Learning to Optimize} (L2O) has emerged as a transformative paradigm that integrates machine learning with numerical optimization to improve computational efficiency, demonstrating particular effectiveness for optimization problems with repetitive structural patterns.
The foundational work of~\cite{gregor2010learning} proposed the learned iterative shrinkage thresholding algorithm, establishing the L2O approach by incorporating trainable parameters into optimization update steps and empirically demonstrating accelerated convergence through data-driven learning. 
Subsequent research has significantly expanded L2O methodologies across various domains, encompassing unconstrained optimization~\cite{andrychowicz2016learning, liu2023towards}, linear optimization~\cite{chen2022representing, li2024pdhg}, combinatorial optimization~\cite{nair2020solving,gasse2022machine}, and constrained nonlinear programs~\cite{donti2020dc3, liang2023low, gaoipm,wu2025constraint, zhang2023jump, liu2023recurrent, liu2021projection, xia2021rnn, du2024long, gao2023novel, liu2023stochastic, xia2022two}. 
They demonstrate remarkable effectiveness when applied to problems with well-defined input distributions. 
Our investigation specifically targets the development of L2O techniques for convex QPs, a crucial subset within constrained optimization frameworks that merits particular attention.

Feasibility and optimality represent two core requirements in constrained L2O frameworks. 
While both supervised and self-supervised learning offer viable approaches, supervised methods critically depend on computationally expensive labeled datasets—a requirement that self-supervised approaches avoid. 
This significant advantage motivates us to focus on self-supervised learning methodologies.
For equality constraints, completion techniques~\cite{donti2020dc3} generate full solutions through algebraic transformations or numerical procedures applied to network outputs~\cite{kim2023self, liang2023low, li2023learning}. 
Inequality constraints are typically addressed via post-processing operations, including gradient-based correction~\cite{donti2020dc3}, Lagrangian dual variable ascent~\cite{kim2023self}, homeomorphic projections~\cite{liang2023low}, and gauge mappings from $\ell_{\infty}$-norm domains~\cite{li2023learning}. 
Another prevalent self-supervised learning strategy incorporates penalty functions into loss formulations to mitigate constraint violation~\cite{donti2020dc3, park2023self, liang2023low}.
For example, the work of~\cite{park2023self} proposed a co-training framework via an augmented Lagrangian method.
While effective at reducing feasibility gaps, penalty methods require meticulous parameter tuning and lack feasibility guarantees.
More importantly, most existing L2O works lack convergence guarantees. 
A notable exception is the work by~\cite{gaoipm}, which integrates approximations based on {Long Short-Term Memory (LSTM)} networks into an inexact IPM framework. 
Their study demonstrates that such learning-based optimization algorithms can attain optimal solutions when the approximation errors remain within acceptable tolerances.

Several pioneering works have explored the integration of L2O techniques with the ADMM. 
The seminal work of~\cite{sun2016deep} first demonstrated this synergy by unrolling ADMM iterations for compressive sensing applications, showing that parameterizing critical operations could significantly improve efficiency for unconstrained inverse problems. 
However, the systematic configuration of convergence-governing parameters in ADMM for convex QPs remains an open research challenge. 
Addressing this gap, the work of~\cite{ichnowski2021accelerating} made substantial advances by developing a reinforcement learning framework to dynamically adapt penalty parameters in OSQP, achieving notable performance improvements through policy gradient optimization.
Complementary approaches have explored warm-start strategies, exemplified by the hybrid architecture in~\cite{sambharya2023end}, which synergizes neural network-generated primal-dual initializations with fixed-iteration Douglas-Rachford splitting. 
These developments collectively highlight the inherent compatibility between ADMM and machine learning paradigms, maintaining theoretical convergence guarantees while exploiting learned priors for computational acceleration.

Although certain L2O algorithms possess theoretical convergence guarantees, a significant theory-practice gap persists in their practical implementation. 
For constrained optimization problems, this discrepancy manifests through substantial objective gaps and measurable constraint violation in real-world applications. 
The L2O paradigm of obtaining approximate solutions represents a departure from exact methods, yet aligns with classical inexact optimization frameworks—where controlled approximations accelerate convergence while maintaining theoretical guarantees, as evidenced in Newton-type methods~\cite{bellavia1998inexact}, IPMs~\cite{eisenstat1994globally}, and ADMM~\cite{bai2025inexact,xie2018inexact}. 
Within ADMM specifically, inexact subproblem solutions have proven particularly effective: The work of~\cite{xie2018inexact} established rigorous relative-error criteria for inexact resolution that enable more aggressive multiplier updates, while~\cite{bai2025inexact} extended these principles to nonconvex problems under proper conditions.
If neural network outputs are designed to satisfy such inexactness criteria, algorithmic convergence can be maintained. 
As demonstrated by~\cite{gaoipm}, empirical adherence to residual bounds in IPMs allows learning-based methods to retain theoretical convergence guarantees.
Building on these insights, we introduce a novel \underline{LSTM}-based framework for generating \underline{i}nexact solutions to \underline{ADMM} subproblems, called \enquote{I-ADMM-LSTM}. 
The proposed architecture employs self-supervised learning on primal-dual residuals across fixed-length iteration windows.
{A safeguard mechanism is incorporated to enforce the prescribed inexactness conditions, thereby guaranteeing convergence.}
The distinct contributions of our work can be summarized as follows:
\begin{itemize}
    {\item \textbf{Learning-based ADMM.} We present a novel neural architecture that seamlessly integrates the traditional ADMM mechanism with deep learning. 
    The framework employs a single LSTM cell—with parameters shared across iterations—to efficiently approximate solutions to the computationally expensive ADMM subproblem. 
    In addition, it adaptively updates key parameters within the inexact ADMM framework using iteration-specific learnable scalars.

    \item \textbf{Safeguard Mechanism}. We establish the conditions that the neural network should satisfy when solving the inexact subproblem to ensure convergence. 
    To strictly enforce these conditions, a safeguard mechanism is introduced: whenever the inexactness criteria are violated, the algorithm reverts to a first-order method to solve the subproblem.
    
    \item \textbf{Theoretical Guarantee.} 
    We establish the convergence of our learning-based inexact ADMM method. 
    In particular, for convex QPs, the algorithm converges to a global optimum, while for nonconvex optimization problems with linear constraints, it converges to a stationary point.}
    
    {\item \textbf{Computational Efficiency.} Extensive evaluations across multiple benchmark problems demonstrate that the proposed method achieves superior optimality and feasibility, along with broader applicability compared to existing L2O approaches.
    Relative to the solvers Gurobi, SCS, and OSQP, our method attains speedups of up to $10\times$, $17\times$, and $14\times$, respectively, while maintaining an objective gap as low as 0.05\%. 
    When extended to nonconvex nonlinear problems with linear constraints, the proposed method achieves speedups of up to $200\times$ over IPOPT, with an objective gap of just 0.5\%.}
\end{itemize}

{The remainder of this paper is structured as follows. Section~\ref{sec:admm_for_convex_qp} presents the theoretical foundations of classical ADMM algorithms for convex QPs and derives convergence conditions for inexact ADMM variants. Section~\ref{sec:our_method} introduces our novel learning-augmented ADMM framework, detailing its architecture and training methodology. 
The theoretical analysis regarding convergence properties of the proposed method appears in Section~\ref{sec:convg_ana}.
Comprehensive numerical experiments on several benchmark problems, demonstrating the practical performance of our approach, are presented in Section~\ref{sec:experiments}. 
Finally, Section~\ref{sec:conclusion} provides concluding remarks that discuss limitations of the proposed method and identify promising directions for future research.}

This paper employs the following mathematical conventions: The symbols $\mathbb{R}$, $\mathbb{R}^n$, and $\mathbb{R}^{n\times m}$ denote the sets of real numbers, $n$-dimensional real vectors, and $n\times m$ real matrices, respectively. 
The cone of $n\times n$ positive semidefinite matrices is represented by $\mathbb{S}^n_{+}$. 
We utilize $\mathbb{I}$ for the indicator function, with $I$ and $0$ designating the identity matrix and zero scalar/vector/matrix of appropriate dimensions. 
For symmetric matrices $A$ and $B$ having matching dimensions, the expressions $A \succ B$ and $A \succeq B$ indicate that $A-B$ is positive definite and positive semidefinite, respectively. 
The standard Euclidean norm on $\mathbb{R}^n$ is written as $\|\cdot\|$. 
For an arbitrary matrix $A$, the notation $\text{Range}(A)$ refers to its range space. 
Finally, given a nonempty closed set $\mathcal{C} \subseteq \mathbb{R}^n$, the Euclidean distance from a point $x$ to $\mathcal{C}$ is expressed as $\text{dist}(x, \mathcal{C}) = \inf_{z \in \mathcal{C}} \|x - z\|$.

\section{ADMM for Convex QPs}\label{sec:admm_for_convex_qp}
This section introduces our convex QP formulation and derives an ADMM framework based on the widely-adopted OSQP solver. 
We subsequently develop an inexact ADMM variant that relaxes the solution accuracy requirements for subproblems, enabling effective incorporation of learning-based optimization approaches.

\subsection{Problem Formulation}

We consider the following convex QP:
\begin{equation}
    \begin{array}{cl}
        \underset{x\in\mathbb{R}^n}{\operatorname{min}} & \frac{1}{2} x^{\top} Q x+p^{\top} x \\
        \text { s.t. } & l \leq A x \leq u,
\end{array}
\label{prob:convex_qp_1}
\end{equation}
where the decision variable $x \in \mathbb{R}^n$ is optimized over a quadratic objective function defined by a positive semidefinite matrix $Q \in \mathbb{S}^n_+$ and a linear term $p \in \mathbb{R}^n$. 
The linear constraints are specified by the matrix $A \in \mathbb{R}^{m \times n}$ with lower and upper bounds $l_i \in \mathbb{R} \cup \{-\infty\}$ and $u_i \in \mathbb{R} \cup \{+\infty\}$, respectively, for each constraint $i = 1, \ldots, m$.
{By introducing an auxiliary variable $z \in \mathbb{R}^m$, the above problem can be equivalently reformulated as
\begin{equation}
    \begin{array}{cl}
\underset{x\in\mathbb{R}^n, z\in\mathbb{R}^m}{\operatorname{min}} & \frac{1}{2} x^{\top} Q x+p^{\top} x \\
\text { s.t. } & A x=z \\
& l \leq z \leq u.
\end{array}
\label{prob:convex_qp_2}
\end{equation}
This formulation is well-suited for solution via ADMM methods and enables efficient projection-based updates for the variable $z$.

Letting $\mathcal{C} = [l, u]$, the optimality conditions for problem \eqref{prob:convex_qp_2} can be expressed as
\begin{align}
&A x = z, \label{eq:opt_conds_1}\\
&Q x + p + A^{T} y = 0, \label{eq:opt_conds_2}\\
&z \in \mathcal{C}, \quad y \in N_{\mathcal{C}}(z),\label{eq:opt_conds_3}
\end{align}
where $y \in \mathbb{R}^m$ denotes the Lagrange multiplier associated with the equality constraint $A x = z$, and $N_{\mathcal{C}}(z)$ represents the normal cone of the set $\mathcal{C}$ at $z$.}

\subsection{OSQP: An ADMM-Based Solver}\label{sec:osqp}
OSQP is a representative ADMM-based solver for convex QPs. 
Given initial primal variables $(x^{0}, z^{0}) \in \mathbb{R}^{n+m}$ and dual variables $y^{0} \in \mathbb{R}^{m}$, and introducing an auxiliary variable $\tilde{z}$ such that $\tilde{z} = z$, the ADMM iterations proceed as follows:

\begin{small}
    \begin{equation}
    \label{eq:admm_step_1}
    \begin{split}
        \left(x^{k+1}, \tilde{z}^{k+1}\right) &\leftarrow \underset{({x}, \tilde{z}): A {x}=\tilde{z}}{\operatorname{argmin}} \; \frac{1}{2} {x}^{\top} Q {x}+p^{\top} {x}+  \frac{\eta}{2} \| {x} -  x^{k} \|_{2}^{2} \\
        &+ \frac{\rho}{2} \left\|\tilde{z}-z^{k}+\rho^{-1} y^{k}\right\|_{2}^{2},
    \end{split}
\end{equation}
\begin{equation}
    \label{eq:admm_step_3}
    z^{k+1} \leftarrow \Pi_{[l, u]}\left(\tilde{z}^{k+1}+\rho^{-1} y^{k}\right), 
\end{equation}
\begin{equation}
    \label{eq:admm_step_4}
    y^{k+1} \leftarrow y^{k}+\rho\left(\tilde{z}^{k+1}-z^{k+1}\right),
\end{equation}
\end{small}where $\eta, \rho>0 $ denote penalty parameters, and $ \Pi_{[l, u]} $ represents the projection operator onto $[l, u]$. 
{Denoting the argument of the projection operator in \eqref{eq:admm_step_3} by $w^k$, the updates $z^{k+1} = \Pi_{[l,u]}(w^k)$ and $y^{k+1} = \rho\big(w^k - \Pi_{[l,u]}(w^k)\big)$ ensure that the optimality condition \eqref{eq:opt_conds_3} is satisfied at each iteration.
This allows us to focus exclusively on the convergence of primal and dual residuals:
\begin{equation}
r_{\text{primal}} \coloneqq A x - z, \quad r_{\text{dual}} \coloneqq Q x + p + A^\top y,
\label{eq:prim_&_dual_residual}
\end{equation}
with OSQP guaranteeing their asymptotic convergence to zero~\cite{banjac2019infeasibility}.}
To address the equality-constrained quadratic subproblem~\eqref{eq:admm_step_1}, we focus on its optimality conditions:
\begin{equation}
    \left[\begin{array}{ccc}
Q+\eta I & 0 & A^{\top} \\
0 & \rho I & -I \\
A & -I & 0 \\
\end{array}\right]\left[\begin{array}{c}
x \\
\tilde{z}\\
v
\end{array}\right]=\left[\begin{array}{c}
-p+\eta x^k \\
\rho z^k - y^k \\
0
\end{array}\right],
\label{eq:sub_opt_conds}
\end{equation}
{where $v$ denotes the dual variable associated with the equality constraint $A x = \tilde{z}$.}
This linear system can be reduced to either: (i) a quasi-definite system solvable by direct methods (e.g., $L^{\top}DL$ factorization), or (ii) a positive definite system amenable to iterative methods (e.g., conjugate gradient).
For direct solution methods, when both $\eta$ and $\rho$ remain fixed, only a single round of matrix factorization is required. 
In practical implementation, while $\eta$ is typically held constant, $\rho$ undergoes limited updates to balance primal-dual residuals; in such cases, only several rounds of matrix factorization are needed. 
This characteristic enables OSQP to achieve superior computational efficiency compared to second-order methods such as IPMs. 
Conversely, indirect methods are better suited for large-scale linear systems and can accommodate frequent updates of $\rho$.

\subsection{Inexact ADMM Framework}
\label{sec:inexact_osqp}

Obtaining an exact solution to subproblem \eqref{eq:admm_step_1} constitutes the primary computational bottleneck in the ADMM iterations \eqref{eq:admm_step_1}–\eqref{eq:admm_step_4}, particularly when dealing with large-scale problems.
To alleviate this issue, we adopt an inexact ADMM approach that computes approximate solutions to~\eqref{eq:admm_step_1}. 
Following the inexact ADMM framework \cite{bai2025inexact}, we reformulate QPs of the form~\eqref{prob:convex_qp_1} as:
\begin{equation}
    \begin{array}{cl}
         \underset{x \in \mathbb{R}^n, z\in\mathbb{R}^m}{\operatorname{min}} & \frac{1}{2} x^{\top} Q x+p^{\top} x + \mathbb{I}_{l\leq z\leq u}(z) \\
         \text{s.t.} & A x=z, \\
\end{array}
\label{prob:convex_qp_3}
\end{equation}
where \(\mathbb{I}_{l\leq z\leq u}(z)\) denotes the indicator function for bound constraints. 
The objective in \eqref{prob:convex_qp_3} consists of a smooth quadratic component \(f(x) \coloneqq \frac{1}{2} x^{\top} Q x + p^{\top} x\) and a nonsmooth term \(g(z) \coloneqq \mathbb{I}_{l \leq z \leq u}(z)\). 
For this two-block convex optimization structure, we employ ADMM with the following iterative steps:
\begingroup
\setlength{\abovedisplayskip}{5pt}
\setlength{\belowdisplayskip}{5pt}
\begin{align}
    z^{k+1} &\leftarrow  \underset{z\in \mathbb{R}^m}{ \text{argmin} } \; \mathcal{L}_{\rho}(x^{k}, z, y^k) , \label{prob:inexact_z_sub_prob} \\
    x^{k+1} &\leftarrow  \underset{x\in \mathbb{R}^n}{ \text{argmin} } \; \mathcal{L}_{\rho}(x, z^{k+1}, y^k)+{\frac{\rho\beta_x}{2}}\|x-x^k\|^2, \label{prob:inexact_x_sub_prob} \\
     y^{k+1} &\leftarrow y^k+\rho\left(Ax^{k+1}-z^{k+1}\right),
     \label{eq:inexact_admm_step_3} 
\end{align}
\endgroup
where $\beta_x>0$ is a penalty parameter, $\mathcal{L}(x, z, y) \coloneqq f(x) + g(z) + y^\top (A x - z)$ and $\mathcal{L}_{\rho}(x, z, y) \coloneqq \mathcal{L}(x, z, y) + \frac{\rho}{2} \|A x - z\|^2$ represent the Lagrangian and augmented Lagrangian functions for Problem~\eqref{prob:convex_qp_3}, respectively.
{The proximal term $\frac{\rho\beta_x}{2}\|x - x^k\|^2$ is introduced in the $x$-update subproblem to facilitate subsequent inexact solves.}

{The $z$-subproblem~\eqref{prob:inexact_z_sub_prob} is a strongly convex QP with an identity quadratic matrix and can therefore be solved analytically via element-wise projection.
In contrast, to reduce the computational cost associated with solving the $x$-subproblem, we propose to solve~\eqref{prob:inexact_x_sub_prob} inexactly, subject to an associated approximation criterion.}
The approximate solution $\tilde{x}^{k+1}$ to the $x$-subproblem \eqref{prob:inexact_x_sub_prob} is required to satisfy the following conditions:
\begin{equation}
    \label{cond:x_subp_cond_1}
    \mathcal{L}_{\rho}(\tilde{x}^{k+1}, z^{k+1}, y^{k}) + \frac{\rho \beta_x}{2} \|\tilde{x}^{k+1} - x^{k}\|^2 \leq \mathcal{L}_{\rho}(x^{k}, z^{k+1}, y^{k}),
\end{equation}
which ensures a monotonic decrease in the objective of \eqref{prob:inexact_x_sub_prob}, together with the gradient condition
{
\begin{equation}
    \label{cond:x_subp_cond_2}
    \|\xi_{x}^{k+1}\|\leq c_{x}\rho \left(\|\tilde{x}^{k+1}-x^k\|+\|z^{k+1}-z^{k}\|\right),
\end{equation}}for some \(c_x > 0\) and \( \xi_x^{k+1} \coloneqq \nabla_x \mathcal{L}_{\rho}(\tilde{x}^{k+1}, z^{k+1}, y^{k}) \).
{This condition ensures that $\tilde{x}^{k+1}$ remains sufficiently close to the exact minimizer.}
{An additional line search procedure is performed after the $y$-update in~\eqref{eq:inexact_admm_step_3}:
\begin{equation}
\label{eq:x_line_search}
    x^{k+1} = \alpha^k \tilde{x}^{k+1} + (1 - \alpha^k)x^k,
\end{equation}}where the step size $\alpha^k \in (0, 2)$ is chosen to satisfy
\begin{equation}
\label{cond:alpha_subp_cond}
    \begin{array}{cl}
         \mathcal{L}_{\rho}(x^{k+1}, z^{k+1}, y^{k+1}) \leq& \mathcal{L}_{\rho}(\tilde{x}^{k+1}, z^{k+1}, y^{k+1}) \\
         & - \delta \rho \|x^{k+1} - \tilde{x}^{k+1}\|^2,
    \end{array}
\end{equation}
with $\delta \in (0, 1)$. Although analogous to the relaxation parameter used in OSQP, $\alpha^k$ here serves the dual role of both accelerating convergence and ensuring algorithmic stability.

Based on the preceding conditions~\eqref{cond:x_subp_cond_1},~\eqref{cond:x_subp_cond_2} and~\eqref{cond:alpha_subp_cond}, we derive an inexact ADMM algorithm (Algorithm~\ref{alg:inexact_admm}) for solving~\eqref{prob:convex_qp_2}. 
At iteration $k$, the algorithm computes the approximate solution $\tilde{x}^{k+1}$ satisfying the specified conditions, followed by a dual variable update and line search procedure to maintain convergence guarantees. 
The iterative process terminates when the composite residual metric
\begin{equation}
    \label{eq:R_k+1}
    R^{k+1} \coloneqq \|\tilde{x}^{k+1} - x^{k}\| + \|z^{k+1} - z^{k}\| + \|A\tilde{x}^{k+1} - z^{k+1}\|
\end{equation}
falls below a predetermined tolerance threshold. 
As proven in Section~\ref{sec:convg_ana}, adherence to these conditions ensures convergence to stationary points of the optimization problem.
\begin{algorithm}[h]  
\caption{Inexact ADMM for convex QPs}  
\label{alg:inexact_admm}  
\textbf{Input}: Initial iterates \(x^0, z^0, y^0\); parameters \(\rho, \beta_x > 0\), \(\alpha \in (0, 2)\)  
\begin{algorithmic}[1]
\STATE \textbf{while}\( R^{k+1} > \epsilon_{\text{tol}} \) \textbf{do}
{
\STATE \hspace{0.5cm}Compute $z^{k+1}$ via~\eqref{prob:inexact_z_sub_prob}\label{step:inexact_admm_step_2}
\STATE \hspace{0.5cm}Identify $\tilde{x}^{k+1}$ such that~(\ref{cond:x_subp_cond_1}) and~(\ref{cond:x_subp_cond_2}) are satisfied\label{step:inexact_admm_step_1}
\STATE \hspace{0.5cm}\(y^{k+1} \leftarrow y^k + \rho(A\tilde{x}^{k+1} - z^{k+1})\)  \label{step:inexact_admm_step_3}
\STATE \hspace{0.5cm}Compute $x^{k+1}$ via~\eqref{eq:x_line_search} such that~\eqref{cond:alpha_subp_cond} is satisfied.
}
\STATE \textbf{end while}
\end{algorithmic}  
\end{algorithm}

\section{Learning-based Inexact ADMM}
\label{sec:our_method}
\begin{figure*}[htbp]
    \centering
    \includegraphics[width=0.93\linewidth]{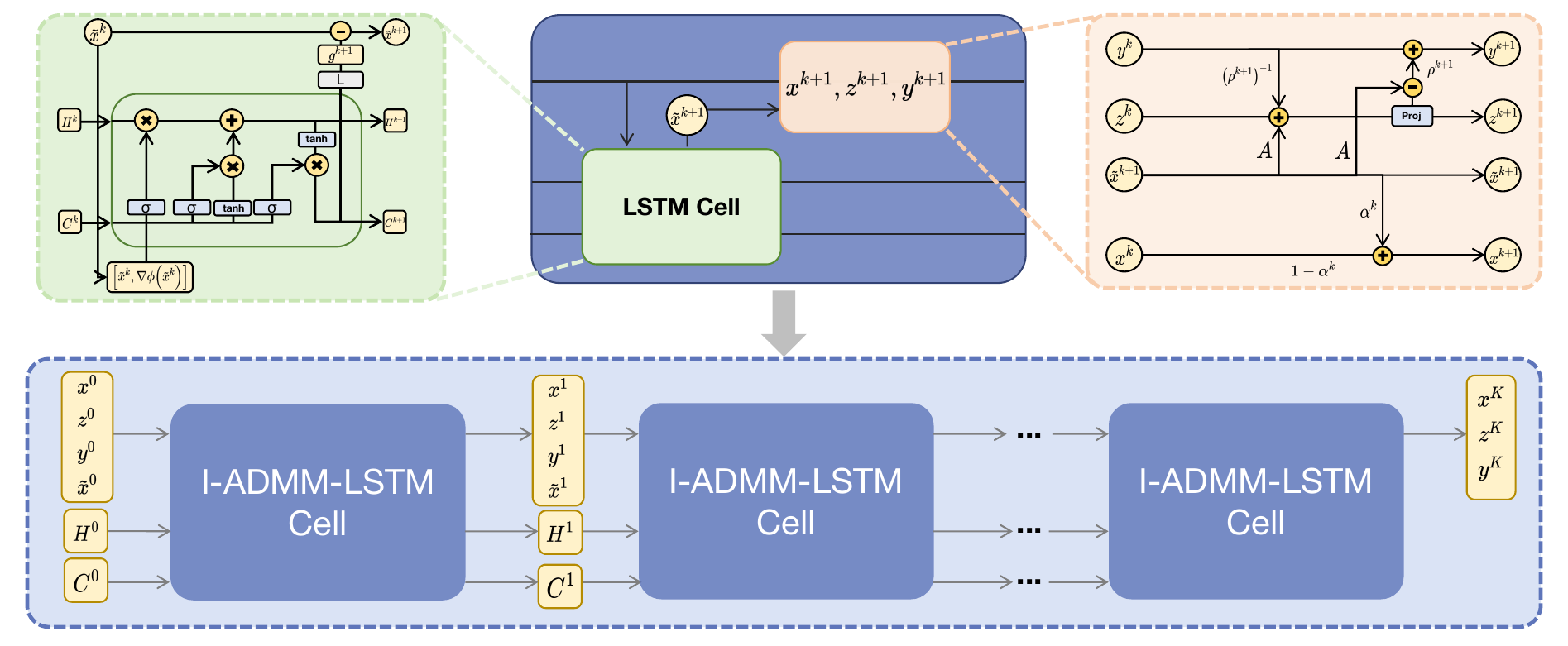}
    \caption{{Illustration of the I-ADMM-LSTM (-Safe) approach. 
    The architecture comprises multiple cells, each receiving and outputting the primal-dual variables $(x, z, y)$, the approximate solution $\tilde{x}$ to subproblem~\eqref{prob:inexact_x_sub_prob}, and the LSTM hidden states $(H, C)$. 
    Each cell contains three components: an LSTM (top-left) that solves~\eqref{prob:inexact_x_sub_prob} to produce $\tilde{x}$; an adaptive ADMM module (top-right) with parameters $\alpha$ and $\rho$ for updating $x$, $y$, and $z$; and an optional safeguard module to ensure strict satisfaction of the inexactness conditions.}}
    \label{fig:I-ADMM-LSTM}
\end{figure*}

Extending the inexact ADMM framework outlined in Section~\ref{sec:admm_for_convex_qp}, we propose a learning-enhanced optimization methodology. 
The core innovation lies in replacing the exact solutions to the ADMM subproblem~\eqref{eq:admm_step_1} with learned approximations generated by LSTM networks, while preserving the overall ADMM iterative structure. 
{To ensure that the inexactness conditions required for convergence are satisfied, an additional safeguard mechanism is introduced.}
Figure~\ref{fig:I-ADMM-LSTM} illustrates the architecture of the proposed method, highlighting the integration of neural approximations with classical optimization components.

\subsection{I-ADMM-LSTM}
{
As discussed in Section~\ref{sec:admm_for_convex_qp}, our primary motivation is to leverage neural networks to provide approximate solutions that avoid explicit matrix factorization, thereby accelerating computation while preserving overall convergence guarantees.
Among all subproblems, only~\eqref{prob:inexact_x_sub_prob} requires matrix factorization when solved exactly. 
An alternative standard approach for solving such unconstrained problems involves first-order methods, particularly gradient-based algorithms. 
Although these methods typically exhibit slow convergence and limited accuracy, they align naturally with our goal of obtaining inexact solutions. Moreover, their iterative structure makes them well-suited for solution trajectory approximation via L2O techniques~\cite{andrychowicz2016learning, gaoipm}, which can significantly improve convergence speed.
}

\begin{figure}[h!]
    \centering
    \includegraphics[width=0.7\linewidth]{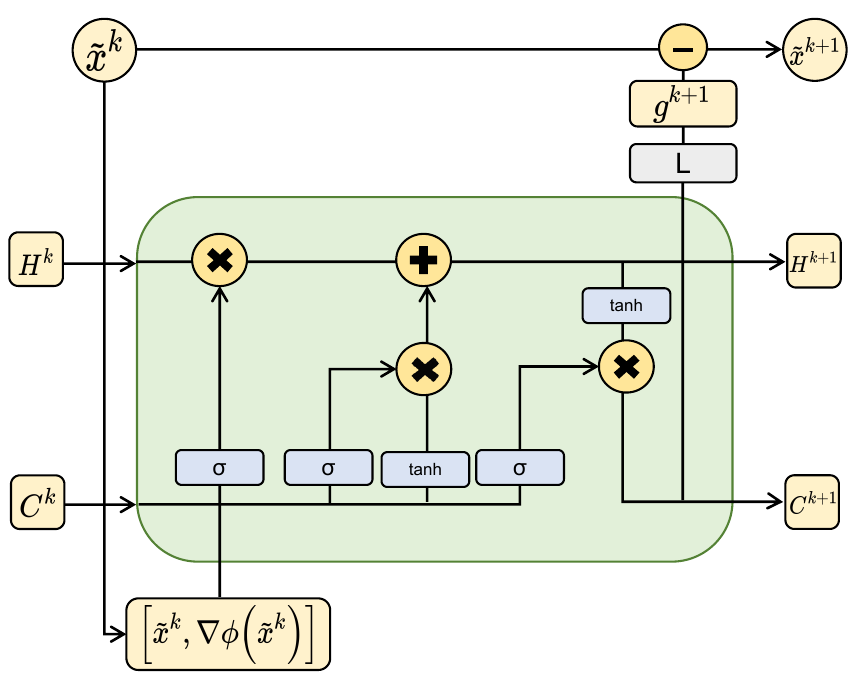}
    \caption{Coordinate-wise LSTM architecture for solving~\eqref{prob:inexact_x_sub_prob}, where \enquote{L} represents the affine transformation layer.}
    \label{fig:lstm}
\end{figure} 

{Our framework employs a single coordinate-wise LSTM cell to generate approximate solutions for subproblem \eqref{prob:inexact_x_sub_prob}, with its architecture detailed in Figure~\ref{fig:lstm}. 
Let the objective function of subproblem \eqref{prob:inexact_x_sub_prob} be defined as
\begin{equation}
\phi(x) := \mathcal{L}_{\rho}(x, z^{k+1}, y^k) + \frac{\rho\beta_x}{2}\|x - x^k\|^2.
\end{equation}
At each iteration, the LSTM cell takes as input the previous approximate solution $\tilde{x}^{k}$ and the gradient $\nabla \phi(\tilde{x}^k)$, and outputs the updated approximation $\tilde{x}^{k+1}$:
\begin{equation}
\label{eq:lstm_update}
\tilde{x}^{k+1} := \operatorname{LSTM}_{\theta}\left(\left[\tilde{x}^{k}, \nabla \phi(\tilde{x}^k)\right]\right).
\end{equation}
Each LSTM cell mimics a first-order algorithm by leveraging gradient information from subproblem \eqref{prob:inexact_x_sub_prob} and accelerates convergence by using the previous approximate solution as an initial guess.}
The parameters of the LSTM cell are shared across all iterations.
By substituting~\eqref{eq:lstm_update} into step~\ref{step:inexact_admm_step_1} of Algorithm~\ref{alg:inexact_admm}, the proposed L2O method can be viewed as a novel recurrent neural network architecture tailored for convex QPs, which we refer to as I-ADMM-LSTM. 
Given initial primal-dual variables $\left( x^{0}, z^{0}, y^{0} \right)$ and LSTM inputs $\tilde{x}^{0}, H^{0}, C^{0}$, the algorithm iterates $K$ times to generate primal-dual solutions $\left(x^{K}, z^{K}, y^{K} \right)$.

Motivated by the convergence criteria of OSQP in~\eqref{eq:prim_&_dual_residual}, we design a self-supervised loss function:
\begin{equation}
    \label{eq:loss_func}
    \min _{\theta} \frac{1}{|\mathcal{B}|} \sum_{B \in \mathcal{B}}\left(\frac{1}{K} \sum_{k=1}^{K} \left\|r_{\text{prim}}^k\right\|+\left\|r_{\text{dual}}^k\right\|\right)_{B},
\end{equation}
where $\theta$ represents trainable network parameters and $\mathcal{B}$ comprises a dataset of identically distributed convex QPs.
This formulation circumvents expensive label generation while ensuring convergence. 
To alleviate memory constraints, we implement truncated backpropagation through time~\cite{liu2023towards}, partitioning the $K$-iteration trajectory into $M$-length subsequences for gradient computation.

\subsection{Adaptive Parameters}
{To enhance algorithmic adaptability, we propose a unified framework for the dynamic adaptation of two key parameters, $\alpha$ and $\rho$, which have a significant impact on the performance of the inexact ADMM algorithm. 
Specifically, the relaxation parameter $\alpha^k \in (0,2)$ must satisfy condition \eqref{cond:alpha_subp_cond} within our inexact ADMM framework. 
Although a line search procedure could be employed to enforce this condition strictly, we opt for a more efficient parameterization using a sigmoid transformation:
\begin{equation}
    \label{eq:adap_alpha}
    \alpha^k = 2 \cdot \sigma({\bar{\alpha}^k}), \quad k=1,\dots,K,
\end{equation}
where $\sigma(\cdot)$ denotes the sigmoid function and $\bar{\alpha}^k \in \mathbb{R}$ is a learnable scalar.

The penalty parameter $\rho > 0$ plays a critical role in determining the overall convergence rate of the algorithm. 
However, the optimal update strategy for $\rho$ remains an open question, and existing approaches primarily rely on heuristic rules. 
In OSQP, $\rho$ is updated only when significant changes are detected to avoid frequent matrix factorizations, which improves computational efficiency but may limit convergence speed~\cite{stellato2020osqp}. 
To strike a balance, we adopt a learnable adaptive scheme:
\begin{equation}
\label{eq:rho_update}
\rho^k := L \cdot \sigma(\bar{\rho}^k),
\end{equation}
where $\bar{\rho}^k \in \mathbb{R}$ is a trainable scalar, and \(L\) denotes the Lipschitz constant of the objective in \eqref{prob:convex_qp_1}, i.e., the maximum eigenvalue of $Q$. 
In practice, $L$ can be estimated by averaging the maximum eigenvalues of $Q$ over the scaled training set.
The inclusion of the scaling factor $L$ is essential for ensuring the theoretical convergence guarantees established in Section~\ref{sec:convg_ana}.

The use of the sigmoid function in both $\alpha^k$ and $\rho^k$ helps mitigate gradient explosion and vanishing issues commonly encountered in deep recurrent neural network architectures.
Crucially, unlike the shared LSTM parameters in Algorithm \ref{alg:inexact_admm}, the parameters $\alpha^k$ and $\rho^k$ are iteration-specific, enabling context-aware adjustments that respond to transient constraint activity. 
This design enhances the algorithm's adaptability across diverse problem geometries while preserving numerical robustness.

To alleviate numerical issues when solving the problem instances, we apply a modified Ruiz scaling procedure prior to optimization, following the approach used in~\cite{stellato2020osqp}. 

}

\subsection{The Safeguard Mechanism}
{

In Section~\ref{sec:inexact_osqp}, the inexact ADMM framework requires that the $x$-subproblem be solved subject to conditions~\eqref{cond:x_subp_cond_1} and~\eqref{cond:x_subp_cond_2}.
Through training, I-ADMM-LSTM effectively minimizes the loss~\eqref{eq:loss_func}, enabling conditions~\eqref{cond:x_subp_cond_1} and~\eqref{cond:x_subp_cond_2} to be approximately satisfied in practice. 
However, strict adherence to these conditions cannot be guaranteed.

To ensure that the approximate solution $\tilde{x}^{k+1}$ strictly satisfies conditions~\eqref{cond:x_subp_cond_1} and~\eqref{cond:x_subp_cond_2}, we introduce an additional safeguard mechanism:
\begin{equation}
\tilde{x}^{k+1} = \left\{
\begin{array}{ll}
\text{Update via}~\eqref{eq:lstm_update}  & \text{if \eqref{cond:x_subp_cond_1}, \eqref{cond:x_subp_cond_2} hold}, \\
\text{Update via}~\operatorname{GD}_{s,T}(\nabla \phi(\tilde{x}^k)) & \text{otherwise},
\end{array}
\right.
\end{equation}
where $\operatorname{GD}_{s,T}$ denotes $T$ steps of gradient descent applied to subproblem \eqref{prob:inexact_x_sub_prob} with step size $s$. 
When the output of the neural network fails to satisfy condition~\eqref{cond:x_subp_cond_1} or~\eqref{cond:x_subp_cond_2}, the algorithm automatically reverts to a conventional optimization method to solve~\eqref{prob:inexact_x_sub_prob}. Since~\eqref{prob:inexact_x_sub_prob} is smooth and strongly convex in $x$, as $T \to \infty$, the gradient descent iterates are guaranteed to satisfy both the sufficient decrease condition~\eqref{cond:x_subp_cond_1} and the gradient condition~\eqref{cond:x_subp_cond_2}~\cite{boyd2004convex, bai2025inexact}.

Moreover, this safeguard is activated only during inference and is not used in training; it therefore does not affect the learning procedure. We refer to this safeguarded variant as I-ADMM-LSTM-Safe, and the complete algorithm is presented in Algorithm~\ref{alg:osqp_lstm}.
}

{
\begin{algorithm}[h]
\caption{I-ADMM-LSTM-Safe for Convex QPs}
\label{alg:osqp_lstm}
\vskip6pt
\textbf{Inputs:} {Initial iterates $x^{0}, z^{0}, y^{0}, \tilde{x}^{0}$, LSTM hidden states $H^{0}, C^{0}$, gradient descent step size $s>0$ and iterations $T$, maximum number of outer iterations $K$}

\textbf{Outputs:} $x^{K}, z^{K}, y^K$
\vskip6pt
\begin{algorithmic}[1]
\FOR{$k \gets 0$ \textbf{to} $K-1$}
{
\STATE $z^{k+1} \leftarrow \Pi_{[l,u]}\left(A\tilde{x}^{k}+({\rho^{k+1}})^{-1} y^{k}\right)$;}
\STATE $\tilde{x}^{k+1}\leftarrow \operatorname{LSTM}_{\theta}\left(\left[\tilde{x}^{k}, \nabla \phi(\tilde{x}^k)\right]\right)$;
{\IF{either \eqref{cond:x_subp_cond_1} or~\eqref{cond:x_subp_cond_2} is not satisfied}
\STATE $\tilde{x}^{k+1}\leftarrow \operatorname{GD}_{s,T}(\nabla \phi(\tilde{x}^k))$
\ENDIF
\STATE $y^{k+1} \leftarrow y^{k}+{\rho^{k+1}}\left(A\tilde{x}^{k+1}-z^{k+1}\right)$;}
\STATE $x^{k+1}\leftarrow {\alpha^{k+1}}\tilde{x}^{k+1}+(1-{\alpha^{k+1}})x^{k}$\label{step:adap_lin_search};
\ENDFOR
\end{algorithmic}
\end{algorithm}
}

\section{Convergence Analysis}
\label{sec:convg_ana}
In this section, we analyze the convergence of Algorithm~\ref{alg:osqp_lstm} under the assumption that I-ADMM-LSTM satisfies conditions~\eqref{cond:x_subp_cond_1},~\eqref{cond:x_subp_cond_2} and~\eqref{cond:alpha_subp_cond} outlined in Section~\ref{sec:inexact_osqp}. 

The quadratic function $f(x) = \frac{1}{2} x^{\top} Q x + p^{\top} x$ has a Lipschitz-continuous gradient, i.e., there exists a constant $L > 0$ such that
\begin{equation}
    \label{eq:lipschitz_ineq}
    \|\nabla f(x_1) - \nabla f(x_2)\| \leq L \|x_1 - x_2\|, \quad \forall x_1, x_2 \in \mathbb{R}^n,
\end{equation}
where $L$ is the Lipschitz constant of $\nabla f$.
We further impose:
\begin{assumption} 
\label{assump:bound}
    For every generated iterate $z^k$, we have $z^k \in \text{Range}(A)$.
\end{assumption}This assumption is readily satisfied, for example, when $A$ has full row rank, in which case $\operatorname{Range}(A)=\mathbb{R}^m$ and any $z^k\in\mathbb{R}^m$ can be represented as $Ax$ for some $x\in\mathbb{R}^n$.
Under this assumption, the dual variable update $y^{k+1} - y^k = \rho \tilde{r}^{k+1}$ lies in $\text{Range}(A)$, yielding:  
\begin{equation}
    \label{eq:dy_ineq}
    \|y^{k+1}-y^{k}\|\leq \left(\sigma^{\text{min}}_{A^{\top}A}\right)^{-\frac{1}{2}}\|A^{\top}\left(y^{k+1}-y^{k}\right)\|
\end{equation}
where $\tilde{r}^{k+1}=A\tilde{x}^{k+1}-z^{k+1}$ and $\sigma_{A^\top A}^{\text{min}}$ is the smallest positive eigenvalue of $A^\top A$. 
To facilitate the analysis, we introduce the following notation: $\tilde{d}_{x}^{k}:=\tilde{x}^{k+1}-x^{k}$,  $\hat{d}_{x}^{k}:=x^{k+1}-\tilde{x}^{k+1}$,${d}_{z}^{k}:=z^{k+1}-z^{k}$ and $d_{y}^{k}:=y^{k+1}-y^{k}$.
We define the potential energy functions as
\begin{equation}
    \label{eq:E_tild}
    \begin{split}
        \tilde{E}^{k+1} &:= \mathcal{L}_{\rho}\left(\tilde{x}^{k+1}, z^{k+1}, y^{k+1}\right)+\tilde{\Gamma}^{k},\\
        {E}^{k+1}&:=\mathcal{L}_{\rho}\left(x^{k+1}, z^{k+1}, y^{k+1}\right)+\Gamma^{k}
    \end{split}
\end{equation}
where
\begin{equation}
    \label{eq:Gamma_tild}
    \tilde{\Gamma}^k := \frac{8\left(1+\tau\right)c^2_x\rho}{ \sigma^{\text{min}}_{A^{\top}A}}\left\|\tilde{d}_{x}^{k}\right\|^{2}+8\left(1+\tau\right)\rho\frac{c_x^2}{\sigma^{\min}_{A^{\top}A}}\|d_z^k\|^2,
\end{equation}
\begin{equation}
    \label{eq:Gamma}
    \Gamma^k := \tilde{\Gamma}^k+\frac{8\left(1+\tau\right)L^2}{\rho \sigma^{\text{min}}_{A^{\top}A}}\left\|\hat{{d}}_{x}^{k}\right\|^{2}
\end{equation}
with $0 < \tau < \delta < 1$.
Then, we can derive the following energy reduction proposition.
\begin{proposition}
    \label{prop:energy_reduc}
    Suppose Assumption~\ref{assump:bound} holds, and conditions~\eqref{cond:x_subp_cond_1}, \eqref{cond:x_subp_cond_2} are satisfied. If the following conditions are satisfied:
    \begin{equation}
    \label{eq:D_x_tild}
    \tilde{\mathcal{D}}_x:=\left[\frac{\beta_x(1-\tau)}{2(1+\tau)}-\frac{2\left(\frac{L}{\rho}+c_x\right)^2+8c_x^2}{\sigma^{\text{min}}_{A^{\top}A}}\right]I \succeq 0,
    \end{equation}
    \begin{equation}
    \label{eq:D_x_hat}
    \hat{\mathcal{D}}_x:=\left[\frac{\delta-\tau}{1+\tau}-\frac{8\left(\frac{L}{\rho}\right)^2}{ \sigma^{\text{min}}_{A^{\top}A}}\right]I \succeq 0,
    \end{equation}
    \begin{equation}
    \label{eq:D_z}
    \mathcal{D}_z:=\left[\frac{(1-\tau)}{2(1+\tau)}-16\frac{c_x^2}{\sigma^{\min}_{A^{\top}A}}\right]I \succeq 0.
    \end{equation}
    then for all $k\geq 1$, we have
    \begin{small}
    \begin{equation}
        \label{eq:E_ineq}
        E^{k+1}\leq E^{k}-\frac{\tau\rho\beta_{x}}{2}\left\|\tilde{d}_{x}^{k}\right\|^{2}-\frac{\tau\rho}{2}\left\|d_{z}^{k}\right\|^{2}-\frac{\tau}{\rho}\left\|d_{y}^{k}\right\|^{2} - \tau\rho\|\hat{d}_x^k\|^2, 
    \end{equation}
    \begin{equation}
        \label{eq:E_tild_ineq}
        \tilde{E}^{k+1} \leq \tilde{E}^k-\frac{\tau\rho\beta_{x}}{2}\left\|\tilde{d}_{x}^{k}\right\|^{2}-\frac{\tau\rho}{2}\left\|d_{z}^{k}\right\|^{2}-\frac{\tau}{\rho}\left\|d_{y}^{k}\right\|^{2}-\tau\rho\|\hat{d}_x^{k-1}\|^2.
    \end{equation}
    \end{small}
\end{proposition}
Readers are referred to Appendix~\ref{app:proofs} for a detailed proof.
The sufficient conditions in~\eqref{eq:D_x_tild}--\eqref{eq:D_z} ensure that the coefficients in the energy descent argument are nonnegative. Given the adaptive update rule for $\rho$ in~\eqref{eq:rho_update}, the practically relevant tunable quantities are mainly $\beta_x$ and $c_x$, while $\tau$ and $\delta$ are auxiliary constants used in the proof. A larger $\beta_x$ strengthens the proximal term in the $x$-subproblem~\eqref{prob:inexact_x_sub_prob}, encouraging the inexact update to remain close to the previous iterate $x^k$. A smaller $c_x$ tightens the relative-error condition in~\eqref{cond:x_subp_cond_2}, requiring the $x$-subproblem to be solved more accurately. Therefore, increasing $\beta_x$ and decreasing $c_x$ generally make these sufficient parameter conditions easier to satisfy, at the cost of more conservative updates or stricter inexactness requirements.
Proposition~\ref{prop:energy_reduc} guarantees the monotonic decrease of the energy sequences $\{E^k\}_{k=1}^{\infty}$ and $\{\tilde{E}^k\}_{k=1}^{\infty}$, forming the foundation for deriving convergence of the proposed algorithm.
A point $(x^*, z^*, y^*)$ is stationary for problem~\eqref{prob:convex_qp_3} if:
\begin{equation}
    \label{cond:optimal_conds}
    0=\nabla f\left(x^{*}\right)+A^{\top} y^{*}, 0 \in \partial g\left(z^{*}\right)-y^{*} \text { and } A x^{*}=z^{*}.
\end{equation}
Then, it is obvious that $(x^k, z^k, y^k)$ is a stationary point if $R^k=0$.
Hence, in the following convergence theorem, we assume $R^k\neq 0$ for all $k$.
\begin{theorem}
    \label{theo:convergence}
    Under Assumption~\ref{assump:bound} and conditions~\eqref{cond:x_subp_cond_1},~\eqref{cond:x_subp_cond_2} and~\eqref{cond:alpha_subp_cond}, if $ \tilde{E}^k$ is bounded below and~\eqref{eq:D_x_tild}~--~\eqref{eq:D_z} hold, there exists $F^*$ such that:
    \begin{small}
        \begin{equation}
        \label{eq:lagran_limit}
        \begin{split}
            \lim _{k \rightarrow \infty} \mathcal{L}\left(x^{k}, z^{k}, y^{k}\right)&=\lim _{k \rightarrow \infty} \mathcal{L}_{\rho}\left(x^{k}, z^{k}, y^{k}\right)=\lim _{k \rightarrow \infty} E^{k}\\
            &=\lim _{k \rightarrow \infty} \tilde{E}^{k}=F^{*},
        \end{split}
    \end{equation}
    \begin{equation}
        \label{eq:paritial_L_limit}
        \begin{split}
            \lim _{k \rightarrow \infty} \operatorname{dist}\left(0, \partial \mathcal{L}\left(x^k, z^k, y^k\right)\right)&=\lim _{k \rightarrow \infty} \operatorname{dist}\left(0, \partial \mathcal{L}_{\rho}\left(x^k, z^k, y^k\right)\right)\\&=0.
        \end{split}
    \end{equation}
    \end{small}
    Let $(x^{*}, z^{*}, y^{*})$ be a limit point of the sequence $\{(x^k,z^k,y^k)\}$. Then $\{(x^k,z^k,y^k)\}$ converges sublinearly to $(x^{*}, z^{*}, y^{*})$, and $(x^{*}, z^{*}, y^{*})$ is an optimal solution to~(\ref{prob:convex_qp_1}).
\end{theorem}
Readers are referred to Appendix~\ref{app:proofs} for a detailed proof. 
Theorem~\ref{theo:convergence} states that Algorithm~\ref{alg:osqp_lstm} achieves sublinear convergence to the optimal solution when mild assumptions are satisfied.

\section{Experiments}
\label{sec:experiments}
We present a systematic evaluation of the proposed \texttt{I-ADMM-LSTM} framework against classical optimization solvers and state-of-the-art L2O methods. 
The comparative analysis focuses on solving convex QPs across multiple benchmarks, quantifying the framework’s performance through objective value gaps, constraint satisfaction residuals, and runtime speedups.

\subsection{Experimental Settings}
\subsubsection{Datasets}
{The convex QPs used in this study are drawn from multiple datasets introduced in~\cite{stellato2020osqp} and~\cite{donti2020dc3}. 
An additional dataset of nonconvex programs from \cite{donti2020dc3} is also included to demonstrate the extensibility of our method to more general nonlinear problems involving linear constraints.
Among these, the datasets denoted as \enquote{Convex QP (RHS)} and \enquote{Nonconvex Program (RHS)} are constructed by perturbing only the right-hand side of the equality constraints.
This setup facilitates a fair comparison with other L2O methods. 
For the remaining datasets, all problem parameters are appropriately perturbed while ensuring that the resulting instances remain feasible.}

\subsubsection{Baseline Algorithms}
In our experiments, we denote our proposed method as \texttt{I-ADMM-LSTM}. The variant incorporating the safeguard mechanism is referred to as \texttt{I-ADMM-LSTM-Safe}. We evaluate its performance against both classical optimization solvers and state-of-the-art L2O algorithms:

\begin{itemize}
    \item \texttt{Gurobi 10.0.3}~\cite{gurobi2023}: A state-of-the-art commercial solver employing advanced numerical techniques for linear, nonlinear, and mixed-integer programming.

    \item \texttt{SCS}~\cite{o2016conic}: An ADMM-based numerical solver for conic programming, which is employed through the CVX framework in this work to streamline implementation.
    
    \item \texttt{OSQP}~\cite{stellato2020osqp}: An ADMM-based solver for convex QPs. 
    To evaluate numerical precision, we compare two configurations: \texttt{OSQP(1E-3)} and \texttt{OSQP(1E-4)}, with absolute/relative tolerances set to $10^{-3}$ and $10^{-4}$, respectively, under a maximum iteration limit of 20,000 and {an adaptive update rule for the penalty parameter $\rho$}.
    
    \item \texttt{DC3}~\cite{donti2020dc3}: A constrained deep learning framework that ensures equality constraints via \emph{completion} steps and enforces inequality constraints through \emph{correction} mechanisms.
    
    \item \texttt{PDL}~\cite{park2023self}: A self-supervised L2O approach that co-trains primal and dual solution networks.
    
    \item \texttt{LOOP-LC}~\cite{li2023learning}: A neural operator leveraging gauge maps to generate feasible solutions for linearly constrained optimization problems.

    \item {\texttt{H-Proj}~\cite{liang2023low}: A method that learns a minimum-distortion homeomorphic mapping from the constraint set to a unit ball using an invertible neural network, followed by a bisection operation to guarantee feasibility of the NN-predicted solution with minimal optimality loss.}

    \item {\texttt{FSNet}~\cite{nguyen2026fsnet}: A L2O method for constrained problems. 
    It first predicts a solution via a multi-layer perceptron, then enforces feasibility by solving an unconstrained reformulation of the constraints using gradient descent.}
\end{itemize}

\subsubsection{Model Settings}
Our model is optimized using the Adam optimizer \cite{kingma2014adam} with a fixed learning rate of $5 \times 10^{-5}$ and a batch size of 4 per task. 
A total of 2,000 perturbed instances are generated and split into 1,880 for training, 20 for validation, and 100 for testing. 
All primal and dual variables, as well as the LSTM state variables in Algorithm~\ref{alg:osqp_lstm}, are initialized to zero.
We employ an early stopping criterion with a patience of 50 epochs, terminating training when no improvement in the objective value is observed for 50 consecutive epochs, provided that both inequality and equality constraint violations remain below predefined tolerance thresholds.
{For \texttt{I-ADMM-LSTM-Safe}, whenever the safeguard mechanism is triggered, we perform gradient descent to solve subproblem \eqref{prob:inexact_x_sub_prob}. The step size and number of iterations are treated as hyperparameters.}

\subsubsection{Configuration.}
All experiments were executed on a computational platform equipped with an NVIDIA RTX A6000 GPU and an Intel Xeon 2.10GHz CPU, utilizing Python 3.10.0 and PyTorch 1.13.1.

\subsection{Computational Results}
\label{sec:comput_res}

We conduct a comprehensive evaluation of \texttt{I-ADMM-LSTM (-Safe)} against baseline methods across multiple datasets. 
All methods are evaluated on 100 test instances, and to ensure a fair comparison with solvers, we report the average performance over sequentially solved instances. 
The evaluation focuses on three key aspects: optimality, feasibility, and efficiency.
The metric \enquote{Obj.} denotes the objective value of the optimization problem, while \enquote{Mean ineq.} and \enquote{Mean eq.} represent the average absolute violations of inequality and equality constraints, respectively. 
For our proposed method, \enquote{Saf.} indicates the average number of times the safeguard mechanism is triggered. 
\enquote{Ite.} reports the average number of iterations for all iterative algorithms, and \enquote{Time (s)} records the sequential solution time in seconds.
Existing L2O methods, due to their limited representational capacity for general nonlinear programming problems, are applicable only to datasets where partial perturbations are applied, such as the “Convex QP (RHS)” dataset. In contrast, our method is capable of handling instances where all parameters are perturbed.
All experimental results are summarized in Table~\ref{tab:synthetic_qp_2}, where the numbers in parentheses under the \enquote{Instance} column denote the number of variables, inequality constraints, equality constraints, lower bounds, and upper bounds, respectively.

Across all datasets, the traditional solvers \texttt{Gurobi}, \texttt{SCS}, and \texttt{OSQP} are all able to guarantee optimality. 
Notably, \texttt{Gurobi} achieves superior computational efficiency due to its commercial-grade implementation, demonstrating average speedups of $3\times$ and $2\times$ over \texttt{SCS} and \texttt{OSQP}, respectively.
Although both \texttt{SCS} and \texttt{OSQP} are based on ADMM, \texttt{SCS} requires more iterations and incurs additional overhead from conic programming reformulations, resulting in inferior performance relative to \texttt{OSQP}.
The relaxed-tolerance variant \texttt{OSQP(1E-3)} exhibits marginally larger objective gaps and feasibility violations compared to \texttt{OSQP(1E-4)}. 
While looser tolerances reduce the number of iterations, they offer limited reduction in matrix factorizations, as $\rho$ updates are triggered only during significant parameter changes. 
Consequently, the overall computation time remains comparable.

On the \enquote{Convex QP (RHS)} dataset, all L2O baseline algorithms returned solutions significantly faster than the solvers. 
However, \texttt{DC3} and \texttt{PDL} exhibited substantial constraint violations across both problem scales. While \texttt{DC3} ensures equality constraint satisfaction via its \enquote{completion} mechanism, its gradient-based \enquote{correction} step fails to fully enforce inequality constraints. 
The end-to-end neural architecture of \texttt{PDL}, lacking explicit feasibility enforcement mechanisms, results in more pronounced violations, particularly for inequality constraints. 
\texttt{LOOP-LC} maintains feasibility through a combination of completion and gauge projection techniques, but at the cost of substantial objective gaps.
{Similarly, \texttt{H-Proj} leverages a homeomorphic projection to effectively guarantee feasibility, yet exhibits considerable objective gaps compared to solvers. 
The recently proposed \texttt{FSNet} achieves strong performance, attaining an objective gap of only $0.7\%$ and maintaining strict feasibility on the order of $10^{-3}$. 
Although its reliance on first-order methods for solving constraint-derived subproblems leads to longer solution times relative to other L2O approaches, it remains significantly faster than traditional solvers.}
Our proposed \texttt{I-ADMM-LSTM} eliminates matrix factorization through learned ADMM parameterization, achieving speedups of $3\times$, $9\times$, and $6\times$ over \texttt{Gurobi}, \texttt{SCS}, and \texttt{OSQP}, respectively, while keeping objective gaps within $1.6\%$ of commercial solvers. 
{Although it exhibits slightly higher equality constraint violations and a marginal degradation in optimality compared to \texttt{FSNet}, it achieves a $33.2\%$ reduction in solution time.
The introduction of the safeguard mechanism further enhances both optimality and feasibility, ultimately enabling \texttt{I-ADMM-LSTM-Safe} to outperform \texttt{FSNet} in solution quality. However, the safeguard is triggered an average of 98.1 times per instance, leading to a noticeable increase in computation time. 
Nevertheless, the safeguarded variant remains faster than \texttt{Gurobi}.}

{Similar performance was observed across the remaining datasets. 
When all parameters are perturbed, the L2O methods are confronted with increased demands in both solution accuracy and adaptability. 
On the \enquote{SVM} dataset, \texttt{I-ADMM-LSTM} achieves an objective gap as low as 0.05\% while maintaining strict feasibility on the order of $10^{-3}$. 
In terms of computational efficiency, \texttt{I-ADMM-LSTM} achieves speedups of $10\times$, $17\times$, and $14\times$ over \texttt{Gurobi}, \texttt{SCS}, and \texttt{OSQP}, respectively, on the \enquote{Portfolio} dataset, while maintaining an objective gap within 1\%.
The proposed safeguard mechanism does not consistently yield improvements. 
For instance, on the \enquote{Random QP} and \enquote{SVM} datasets, a slight degradation in optimality is observed. 
This can be attributed to two factors: first, the limited gradient descent procedure sometimes fails to deliver sufficiently accurate solutions; second, the iterative trajectory of the learning-based inexact ADMM, optimized via a loss function, does not always align with that of conventional optimization methods.}

\begin{table*}[!t]
	\centering
	\caption{Comparative results on synthetic datasets.}
	\label{tab:synthetic_qp_2}
        \resizebox{0.9\textwidth}{!}{
		\begin{tabular}{clcccccc}
		\toprule
		Instance & Method & Obj.  & Mean ineq. $\downarrow$ &Mean eq. $\downarrow$ & \# Saf. & Ite.  & Time (s) $\downarrow$  \\
        \midrule
        \multirow{11}{*}{\makecell{Convex QP (RHS) \\ (1,000, 500, 500, 0, 0)}}& \texttt{Gurobi} & -152.701 & 0.000 & 0.000 & - & 14.3 & 0.824  \\
        &\texttt{SCS} & -152.701 & 0.000 & 0.000 & - & 75.0 & 2.285  \\
		&\texttt{OSQP(1E-4)} & -152.701 & 0.000 & 0.000 & - & 33.6 & 1.491  \\
		&\texttt{OSQP(1E-3)} &  -152.705 &  0.002 & 0.001 & - & 28.3 & 1.456 \\
        & \texttt{DC3} & -112.894  &  0.175 & 0.000 & - & - & 0.056   \\
        & \texttt{PDL} & -131.568  &  0.002 & 0.555 & - & - & 0.025  \\
        & \texttt{LOOP-LC} & -108.525  &  0.000 & 0.000 & - & - & 0.042  \\
        & {\texttt{H-Proj}} & {-106.525} & {0.000} & {0.000} & {-} & {-} & {0.084} \\
        & {\texttt{FSNet}} & {-151.573} & {0.000} & {0.000} & {-} & {-} & {0.397} \\
        &{\texttt{I-ADMM-LSTM}} & {-150.265}  & {0.000}  & {0.002} & {0.0} & {150.0} &  {0.265} \\
        &{\texttt{I-ADMM-LSTM-Safe}} & {-151.701}  & {0.000}  & {0.001} & {98.1} & {150.0} &  {0.705} \\
            \midrule
            \multirow{6}{*}{\makecell{Convex QP (ALL) \\ (1,000, 500, 500, 0, 0)}}& \texttt{Gurobi} & -164.035 & 0.000 & 0.000 & - & 14.3 & 0.831 \\
            & \texttt{SCS} & -166.014 & 0.000 & 0.000 & - & 75.0 & 2.353 \\
		&\texttt{OSQP(1E-4)} & -164.035 & 0.000 & 0.000 & - & 35.9 & 1.558  \\
		&\texttt{OSQP(1E-3)} & -164.039 &  0.002 & 0.001 & - & 30.4 & 1.504  \\
        &{\texttt{I-ADMM-LSTM}} & {-161.010}  &  {0.001} & {0.008} & {0.0} & {150.0} & {0.255}  \\
		&{\texttt{I-ADMM-LSTM-Safe}} & {-162.820} &  {0.000} & {0.004} & {101.1} & {150.0} & {0.712}  \\
        \midrule
            \multirow{6}{*}{\makecell{Random QP \\ (1,000, 2,000, 0, 0, 0)}}& \texttt{Gurobi} & -2.714 & 0.000 & - & - &  12.9 & 3.934 \\
            & \texttt{SCS} & -2.714 & 0.000 & - & - & 191.7 & 8.060 \\
		&\texttt{OSQP(1E-4)} & -2.714 & 0.000 & - & - & 93.9 & 10.765 \\
		&\texttt{OSQP(1E-3)} & -2.714 & 0.000 & - & - & 79.4 & 9.968 \\ 
        & {\texttt{I-ADMM-LSTM}} & {-2.630} & {0.000} & {-} & {0.0} & {600.0} & {0.975} \\
        & {\texttt{I-ADMM-LSTM-Safe}} & {-2.592} & {0.000} & {-} & {392.3} & {600.0} & {2.667} \\
            \midrule
            \multirow{6}{*}{\makecell{Equality QP \\ (1,000, 0, 500, 0, 0)}}& \texttt{Gurobi} & 249.722 & - & 0.000 & - & 5.3 & 0.794 \\
            & \texttt{SCS} & 249.722 & - & 0.000 & - & 25.0 & 1.509 \\
		&\texttt{OSQP(1E-4)} & 249.732 & - & 0.000 & - & 36.1 & 1.245 \\
		&\texttt{OSQP(1E-3)} & 249.740 & - & 0.000 & - & 31.7 & 1.234 \\
        &{\texttt{I-ADMM-LSTM}} &  {250.237} & {-} & {0.001} & {0.0} & {150.0} & {0.264} \\
		&{\texttt{I-ADMM-LSTM-Safe}} &  {249.749} & {-} & {0.000} & {34.9} & {150.0} & {0.578} \\
        \midrule
            \multirow{6}{*}{\makecell{{Portfolio} \\ (1,000, 0, 1, 1,000, 1,000)}}& {\texttt{Gurobi}} & {-0.702} & {0.000} & {0.000} & {-} & {14.7} & {0.909} \\
            & {\texttt{SCS}} & {-0.702} & {0.000} & {0.000} & {-} & {115.8} & {1.621} \\
            & {\texttt{OSQP(1E-4)}} & {{-0.716}} & {0.000} & {0.000} & {-} & {72.2} & {1.334} \\
		&{\texttt{OSQP(1E-3)}} & {-0.757} & {0.000} & {0.001} & {-} & {33.0} & {0.724} \\
        &{\texttt{I-ADMM-LSTM}} & {-0.710}  & {0.000} & {0.000} & {0.0} & {50.0} & {0.095} \\
		&{\texttt{I-ADMM-LSTM-Safe}} & {-0.702}  & {0.000} & {0.000} & {40.9} & {50.0} & {0.620} \\
            \midrule
            \multirow{6}{*}{\makecell{SVM \\ (1,500, 500, 0, 0, 0)}}& \texttt{Gurobi} & 547.880 & 0.000 & - & - & 6.8 & 0.236 \\
            & \texttt{SCS} & 547.880 & 0.000 & - & - & 50.0 & 2.176 \\
            &\texttt{OSQP(1E-4)} & 547.859 & 0.000 & - & - & 51.0 & 0.227 \\
		&\texttt{OSQP(1E-3)} & 548.226 & 0.000 & - & - & 36.7 & 0.171 \\
        &{\texttt{I-ADMM-LSTM}} &  {547.593} & {0.000} & {-} & {0.0} & {50.0} & {0.119} \\
		&{\texttt{I-ADMM-LSTM-Safe}} &  {549.091} & {0.000} & {-} & {18.7} & {50.0} & {0.449} \\
		\bottomrule
	\end{tabular}}
\end{table*}

\subsection{Performance Analysis}
\label{sec:perform_analysis}
\subsubsection{Convergence Profiling}
\begin{figure*}[!t]
    \centering
\includegraphics[width=0.9\linewidth]{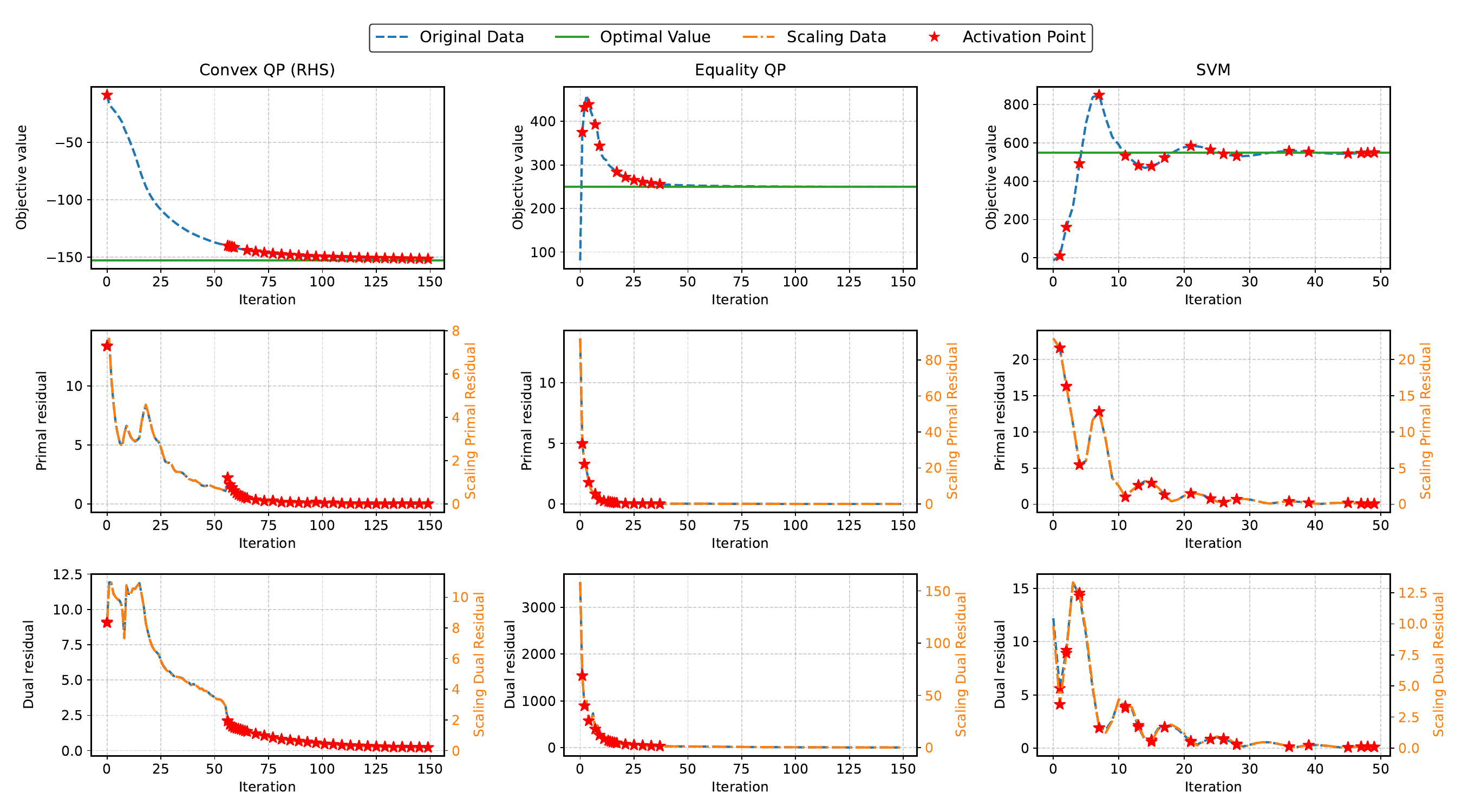}
    \caption{{Convergence behavior of \texttt{I-ADMM-LSTM-Safe} on representative problem instances, showing the evolution of the objective value as well as primal and dual residuals. 
    The blue curves correspond to the raw output of \texttt{I-ADMM-LSTM-Safe}. 
    In the residual plots, the orange lines indicate the residuals after applying Ruiz scaling, while the red markers highlight iterations where the safeguard mechanism is activated.}}
    \label{fig:obj_prim_dual_res}
\end{figure*}

{
To comprehensively evaluate the convergence behavior of our algorithm, we focus on the safeguarded variant, \texttt{I-ADMM-LSTM-Safe}, and analyze its performance on three representative datasets: \enquote{Convex QP (RHS)}, \enquote{Equality QP}, and \enquote{SVM}. Specifically, we examine the evolution of the objective value, as well as the primal and dual residuals, with results summarized in Figure~\ref{fig:obj_prim_dual_res}.

Across all three datasets, \texttt{I-ADMM-LSTM-Safe} exhibits consistent convergence, with the objective values steadily approaching the optimum while both primal and dual residuals converge toward zero.
The orange curves in the second and third rows depict the residuals after applying Ruiz scaling. This scaling effectively mitigates numerical ill-conditioning, particularly for the \enquote{Equality QP} dataset, where the dual residual range is reduced from $[0, 3500]$ to $[0, 150]$, bringing it to a similar scale as the primal residual $[0,100]$. 
Our loss function~\eqref{eq:loss_func} incorporates both primal and dual residuals, which encourages the model to balance these two components during training irrespective of their original magnitudes.

Red stars in the figure mark the iterations where the safeguard mechanism is activated; overlapping markers have been thinned for visual clarity. 
The activation patterns vary across datasets, reflecting the distinct solution characteristics of each problem. 
For the \enquote{Convex QP (RHS)} dataset, safeguard activations are concentrated in the later stages of convergence, indicating that the algorithm encounters difficulties satisfying the inexactness conditions near the optimum. 
In the \enquote{Equality QP} case, activations occur primarily during phases of rapid convergence, suggesting that while the algorithm makes swift progress, it may temporarily violate the required conditions. 
For the \enquote{SVM} dataset, safeguard activations appear at each turning point of the convergence trajectory, where fluctuations in the solution path are detected and subsequently corrected by the safeguard.
Overall, these observations suggest that the safeguard mechanism is more likely to be activated near high-accuracy regimes, during rapid convergence phases, or around turning points of the learned trajectory, where small prediction errors may prevent the LSTM output from providing a sufficiently stable inexact ADMM step.
}

{


\subsubsection{Runtime Analysis}

{
To analyze the computational cost of our algorithm, we decompose the runtime of \texttt{I-ADMM-LSTM-Safe} into four components: (i) \enquote{Scaling}: the time for Ruiz scaling and solution recovery; (ii) \enquote{LSTM}: the inference time of the LSTM network; (iii) \enquote{SF-Cond}: the time required to check whether $\tilde{x}^k$ satisfies the inexactness conditions \eqref{cond:x_subp_cond_1} and \eqref{cond:x_subp_cond_2}; and (iv) \enquote{SF-GD}: the time spent on gradient descent when the safeguard is triggered. 
For \texttt{I-ADMM-LSTM}, only the first two components are relevant. 
We compare CPU and GPU runtimes on the \enquote{Convex QP (ALL)} and \enquote{Portfolio} datasets, with results shown in Figure~\ref{fig:cpu_gpu_comp}. 
As expected, learning-based methods benefit significantly from GPU parallelism due to their reliance on matrix operations. 
For instance, the LSTM runtime increases from $0.238$~s and $0.074$~s on GPU to $6.588$~s and $1.517$~s on CPU for the two datasets, respectively. 
In contrast, the safeguard mechanism, which is dominated by iterative gradient descent, exhibits only a modest increase of around 20\% on CPU. 
The scaling step also shows a notable increase on CPU, as it involves substantial matrix computations.
}

\begin{figure}[!t]
    \centering
    \includegraphics[width=1.0\linewidth]{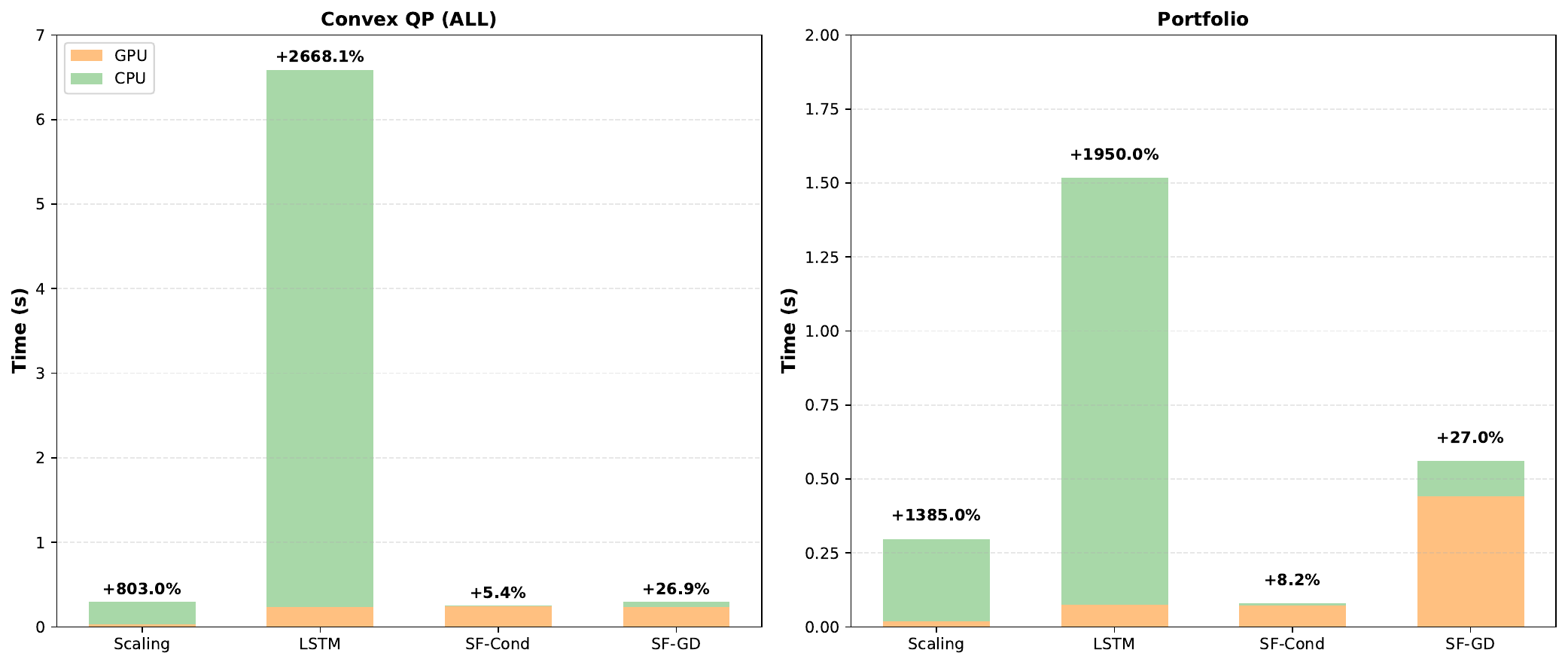}
    \caption{{Runtime breakdown of the proposed algorithm on CPU and GPU for two representative datasets: (a) \enquote{Convex QP (ALL)} and (b) \enquote{Portfolio}. The percentages above the CPU bars indicate the runtime increase relative to GPU execution.}}
    \label{fig:cpu_gpu_comp}
\end{figure}

{
To ensure a fair comparison with traditional solvers, we report in Table~\ref{tab:synthetic_qp_2} the sequential solution time (with batch size $=1$) for all methods. 
However, learning-based approaches, including ours, are inherently capable of solving multiple optimization instances in parallel. 
When executed on a GPU, this parallelization yields significant computational advantages over sequential solvers, as also noted in \cite{nguyen2026fsnet}.
We further exploit the parallel processing capabilities of \texttt{I-ADMM-LSTM (-Safe)} and analyze how its solution time scales with increasing batch size. 
The results are presented in Figure~\ref{fig:para_time}. 
Across all six convex QP datasets, the runtime of \texttt{I-ADMM-LSTM (-Safe)} decreases substantially as the batch size increases. 
This downward trend begins to moderate as the batch size approaches 32, due to the progressive saturation of GPU computational capacity.
}

\begin{figure*}[!t]
    \centering
\includegraphics[width=0.8\linewidth]{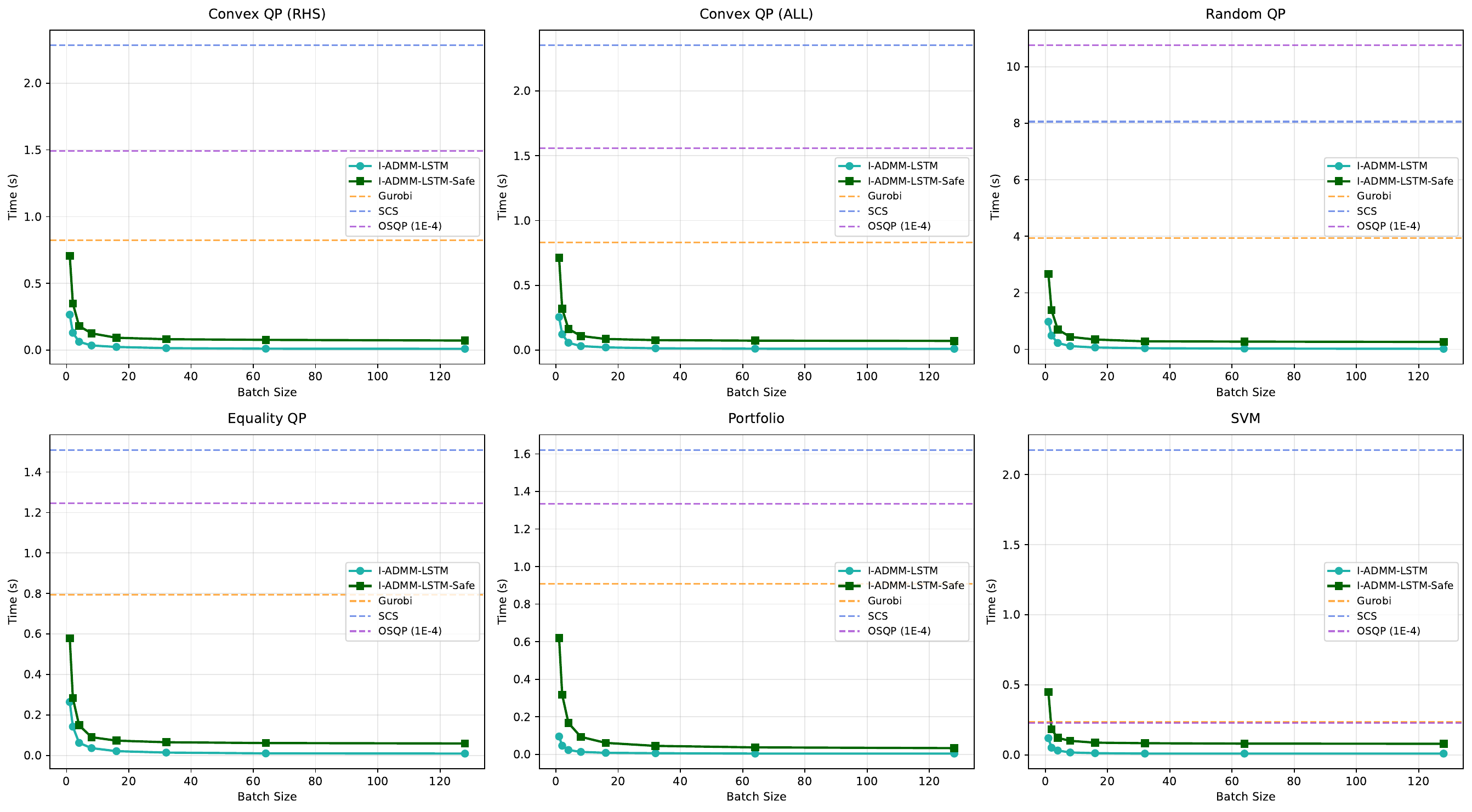}
    \caption{{The time required for \texttt{I-ADMM-LSTM (-Safe)} to solve multiple optimization problems in parallel.}}
    \label{fig:para_time}
\end{figure*}


\subsection{Extension to Nonconvex Problems with Linear Constraints}

\begin{table*}[!t]
	\centering
	\caption{Comparative performance evaluation on \enquote{Nonconvex Programs (RHS)} and \enquote{Nonconvex Programs (ALL)}.}
	\label{tab:synthetic_nonconvex_program}
        \resizebox{0.95\textwidth}{!}{
		\begin{tabular}{clcccccc}
		\toprule
		Instance & Method & Obj.  & Mean ineq. $\downarrow$ &Mean eq. $\downarrow$ & \# Saf. &  Ite. & Time (s) $\downarrow$  \\
            \midrule
            \multirow{7}{*}{\makecell{Nonconvex Program (RHS) \\ (1,000, 500, 500)}}
            &{\texttt{IPOPT}} & {-125.132} & {0.000} & {0.000} & {-} & {14.9} & {37.778} \\
        & {\texttt{DC3}} & {-113.381} & {0.495}  & {0.000} & {-} & {-} & {0.048}  \\
        & {\texttt{PDL}} & {-73.436} & {0.000} & {0.504} & {-} & {-} & {0.016} \\
        & {\texttt{LOOP-LC}} & {109.433} & {0.000} & {0.000} & {-} & {-} & {0.028}  \\
        & {\texttt{FSNet}} & {-123.980} & {0.000} & {0.000} & {-} & {-} & {0.419} \\
        &{\texttt{I-ADMM-LSTM}} & {-124.385} & {0.000}  & {0.001}  & {0.0} & {150.0} & {0.260}  \\
        &{\texttt{I-ADMM-LSTM-Safe}} & {-124.767} & {0.000}  & {0.001}  & {82.6} & {150.0} & {0.935}  \\
            \midrule
            \multirow{3}{*}{\makecell{Nonconvex Program (ALL) \\ (1,000, 500, 500)}}& {\texttt{IPOPT}} & {-127.491} & {0.000} & {0.000} & {-} & {17.4} & {50.830}  \\
            &{\texttt{I-ADMM-LSTM}} & {-124.736} & {0.001} & {0.007} & {0.0} & {150.0} & {0.251}   \\
            &{\texttt{I-ADMM-LSTM-Safe}} & {-125.502} & {0.000} & {0.002} & {59.2} & {150.0} & {0.726}   \\
            
		\bottomrule
	\end{tabular}}
\end{table*}

Although \texttt{I-ADMM-LSTM (-Safe)} is designed to find global optima for convex QPs, it can also be applied to other optimization problems with smooth objective functions—including nonconvex ones—subject to linear constraints. 
Since only the form of the objective function changes, as long as conditions \eqref{cond:x_subp_cond_1} and \eqref{cond:x_subp_cond_2} are satisfied, and the additional line search satisfies condition \eqref{cond:alpha_subp_cond}, our theoretical results remain valid and guarantee convergence to a stationary point.

In this section, we consider a nonconvex optimization problem with 1,000 variables, 500 inequality constraints, and 500 equality constraints, as implemented in~\cite{donti2020dc3,park2023self,li2023learning}:
\begin{equation}
\label{prob:ncvx_program}
\begin{array}{cl}
\underset{x \in \mathbb{R}^{n} }{ \min}   & \frac{1}{2} x^{\top} Q x + p^{\top} \sin(x) \\
\text { s.t. } & l \leq A x \leq u.
\end{array}
\end{equation}
where $\sin(\cdot)$ is applied element-wise to the vector $x$.
We compare our method with \texttt{DC3}, \texttt{PDL}, and \texttt{LOOP}. 
Since the problem is nonconvex and nonlinear, we also compare against IPOPT as a strong local optimization baseline. 
The results are summarized in Table~\ref{tab:synthetic_nonconvex_program}.
Similar to the observations on the \enquote{Convex QP (RHS/ALL)} datasets, all L2O baselines exhibit comparable optimality and feasibility performance but are limited to settings with partial parameter perturbations. In contrast, \texttt{I-ADMM-LSTM (-Safe)} is capable of handling a wide range of perturbed problem instances while maintaining strong feasibility and competitive local solution quality, achieving an objective gap as low as 0.5\% relative to IPOPT. Moreover, given that nonconvex nonlinear problems are inherently more challenging for traditional solvers, \texttt{I-ADMM-LSTM} achieves speedups of up to nearly \(200\times\) over IPOPT while maintaining this small objective gap.

}

\FloatBarrier
\section{Conclusions}\label{sec:conclusion}
{In this study, we propose I-ADMM-LSTM, a learning-based inexact ADMM framework that approximates computationally expensive ADMM subproblem solutions using a specialized LSTM network architecture. 
We establish the inexactness conditions that the approximate solutions must satisfy and introduce a safeguard mechanism to enforce them, thereby ensuring the convergence of the proposed L2O method.
Extensive computational experiments on multiple convex QPs validate the efficiency of I-ADMM-LSTM in generating high-quality globally optimal solutions. 
Similar performance is observed when the method is extended to nonconvex and nonlinear problems with linear constraints, where it aims for stationary solutions and empirically attains competitive local solution quality while achieving faster solve times compared to traditional solvers.
Future work will aim to extend the proposed framework to problems with nonlinear constraints, while preserving both solution quality and convergence guarantees.}



\appendices

\section{Proofs}
\label{app:proofs}
Before proving Proposition~\ref{prop:energy_reduc}, we first establish the following auxiliary lemma.
\begin{lemma}
    \label{lemma:Ady}
    Under Assumption~\ref{assump:bound} and Condition~\eqref{cond:x_subp_cond_2}, the following inequality holds for all $k \geq 1$: 
    \begin{equation}
        \label{eq:lemma_1}
        \begin{array}{cl}
            \|A^{\top}d_y^k\|^2 \leq& 2\left(L+c_x\rho\right)^2\|\tilde{d}_{x}^{k}\|^2+8L^2\|\hat{d}_{x}^{k-1}\|^2\\
            &+8c_x^2\rho^2\left(\|\tilde{d}_x^{k-1}\|^2+\|d_z^k\|^2+\|d_z^{k-1}\|^2\right)
        \end{array}
    \end{equation}
\end{lemma}
\subsection{Proof of Lemma~\ref{lemma:Ady}}
\label{app:proof_lemma}
By the definition of the gradient $\xi_{x}^{k+1}=\nabla_{x} \mathcal{L}_{\rho}\left(\tilde{x}^{k+1}, z^{k+1}, y^{k}\right)$, we have
\begin{equation}
    \label{eq:xi}
    \begin{array}{cll}
        \xi_{x}^{k+1} & = &\nabla f\left(\tilde{x}^{k+1}\right)+A^{\top}y^{k} +\rho A^{\top}\left(A\tilde{x}^{k+1} - z^{k+1}\right)\\
         & =&\nabla f\left(\tilde{x}^{k+1}\right)+A^{\top}y^{k}+\rho A^{\top}\tilde{r}^{k+1}\\
    \end{array}
\end{equation}
where $\tilde{r}^{k+1}=A\tilde{x}^{k+1} - z^{k+1}$. It follows that
\begin{equation}
    \label{eq:Ay}
    A^{\top} y^{k}=\xi_{x}^{k+1} - \nabla f\left(\tilde{x}^{k+1}\right)- \rho A^{\top} \tilde{r}^{k+1}
\end{equation}
Substituting the dual update $y^{k+1} = y^k + \rho \tilde{r}^{k+1}$ into \eqref{eq:Ay} gives
\begin{equation}
    A^{\top} y^{k+1}=\xi_{x}^{k+1} - \nabla f\left(\tilde{x}^{k+1}\right)
\end{equation}
Letting $d_y^k=y^{k+1}-y^k$, we obtain
\begin{equation}
    \begin{array}{cl}
       A^{\top} d^{k}_y = &\xi_{x}^{k+1} - \xi_{x}^{k} +\nabla f\left(\tilde{x}^{k}\right) - \nabla f\left(\tilde{x}^{k+1}\right)
    \end{array}
\end{equation}
Taking norms on both sides and applying the Cauchy--Schwarz inequality yield
\begin{equation}
    \begin{aligned}
\left\|A^{\top} d_{y}^{k}\right\|^{2} & = \|\nabla f\left(\tilde{x}^{k}\right) - \nabla f\left(\tilde{x}^{k+1}\right) + \xi_{x}^{k+1} - \xi_{x}^{k}\|^2 \\
& \leq \left[L\|\tilde{d}_{x}^{k}+\hat{d}_{x}^{k-1}\|+c_x\rho\left(\|\tilde{d}_x^k\|+\|\tilde{d}_x^{k-1}\|\right.\right.\\
&+\left.\left.\|d_z^k\|+\|d_z^{k-1}\|\right) \right]^2\\
& \leq \left[\left(L+c_x\rho\right)\|\tilde{d}_{x}^{k}\|+L\|\hat{d}_{x}^{k-1}\|\right.\\
&+\left.c_x\rho\left(\|\tilde{d}_x^{k-1}\|+\|d_z^k\|+\|d_z^{k-1}\|\right) \right]^2\\
& \leq 2\left(L+c_x\rho\right)^2\|\tilde{d}_{x}^{k}\|^2+8L^2\|\hat{d}_{x}^{k-1}\|^2\\
&+8c_x^2\rho^2\left(\|\tilde{d}_x^{k-1}\|^2+\|d_z^k\|^2+\|d_z^{k-1}\|^2\right) \\
\end{aligned}
\end{equation}

\subsection{Proof of Proposition~\ref{prop:energy_reduc}}
\label{app:proof_proposition}
We now prove Proposition~\ref{prop:energy_reduc}. For fixed \(x^k\) and \(y^k\), the function \(\mathcal{L}_{\rho}(x^k, z, y^k)\) in subproblem \eqref{prob:inexact_z_sub_prob} is \(\rho\)-strongly convex in \(z\). Therefore,
\begin{equation}
\mathcal{L}_{\rho}\left(x^k, z^{k+1}, y^{k}\right) + \frac{\rho}{2}\left\|z^{k+1} - z^{k}\right\|^{2} \leq \mathcal{L}_{\rho}\left(x^k, z^{k}, y^{k}\right).
\label{eq:strong_convex_z}
\end{equation}
Combining \eqref{eq:strong_convex_z} with condition \eqref{cond:x_subp_cond_1} and the dual residual bound \eqref{eq:dy_ineq}, we obtain the following decomposition of the augmented Lagrangian difference:
\begin{equation}
    \label{eq:theo_1_proof_1}
    \begin{aligned}
& \mathcal{L}_{\rho}\left(\tilde{x}^{k+1}, z^{k+1}, y^{k+1}\right)-\mathcal{L}_{\rho}\left(x^{k}, z^{k}, y^{k}\right) \\
= & \mathcal{L}_{\rho}\left(\tilde{x}^{k+1}, z^{k+1}, y^{k+1}\right)-\mathcal{L}_{\rho}\left(\tilde{x}^{k+1}, z^{k+1}, y^{k}\right)\\
&+\mathcal{L}_{\rho}\left(\tilde{x}^{k+1}, z^{k+1}, y^{k}\right)  -\mathcal{L}_{\rho}\left(x^{k}, z^{k+1}, y^{k}\right)\\
&+\mathcal{L}_{\rho}\left(x^{k}, z^{k+1}, y^{k}\right)-\mathcal{L}_{\rho}\left(x^{k}, z^{k}, y^{k}\right) \\
\leq & \frac{1+\tau}{\rho}\left\|d_{y}^{k}\right\|^{2}-\frac{\rho}{2}\left\|d_{z}^{k}\right\|^{2}-\frac{\rho\beta_{x}}{2}\left\|\tilde{d}_{x}^{k}\right\|^{2}-\frac{\tau}{\rho}\left\|d_{y}^{k}\right\|^{2} \\
\leq & \frac{1+\tau}{\rho \sigma^{\text{min}}_{A^{\top}A}}\left\|A^{\top} d_{y}^{k}\right\|^{2}-\frac{\rho}{2}\left\|d_{z}^{k}\right\|^{2}-\frac{\rho\beta_{x}}{2}\left\|\tilde{d}_{x}^{k}\right\|^{2}-\frac{\tau}{\rho}\left\|d_{y}^{k}\right\|^{2}.
\end{aligned}
\end{equation}
Moreover, Lemma~\ref{lemma:Ady} implies
\begin{equation}
\label{eq:dy_bound_for_energy}
\begin{aligned}
& \frac{1+\tau}{\rho \sigma^{\text{min}}_{A^{\top}A}}\left\|A^{\top} d_{y}^{k}\right\|^{2} \\
\leq & \frac{1+\tau}{\rho \sigma^{\text{min}}_{A^{\top}A}}\left[2\left(L+c_x\rho\right)^2\|\tilde{d}_{x}^{k}\|^2+8L^2\|\hat{d}_{x}^{k-1}\|^2 \right.\\
&\left.+ 8c_x^2\rho^2\left(\|\tilde{d}_x^{k-1}\|^2+\|d_z^k\|^2+\|d_z^{k-1}\|^2\right)\right].
\end{aligned}
\end{equation}
Substituting \eqref{eq:dy_bound_for_energy} into \eqref{eq:theo_1_proof_1} and using \eqref{cond:alpha_subp_cond}, we obtain
\begin{equation}
\label{eq:theo_1_proof_3}
\begin{aligned}
& \mathcal{L}_{\rho}\left(x^{k+1}, z^{k+1}, y^{k+1}\right)-\mathcal{L}_{\rho}\left(x^{k}, z^{k}, y^{k}\right) \\
\leq & \mathcal{L}_{\rho}\left(\tilde{x}^{k+1}, z^{k+1}, y^{k+1}\right)-\mathcal{L}_{\rho}\left(x^{k}, z^{k}, y^{k}\right) - \delta \rho\|x^{k+1}-\tilde{x}^{k+1}\|^2 \\
\leq & \frac{1+\tau}{\rho \sigma^{\text{min}}_{A^{\top}A}}\left[2\left(L+c_x\rho\right)^2\|\tilde{d}_{x}^{k}\|^2+8L^2\|\hat{d}_{x}^{k-1}\|^2\right.\\
&\left.+8c_x^2\rho^2\left(\|\tilde{d}_x^{k-1}\|^2+\|d_z^k\|^2+\|d_z^{k-1}\|^2\right)\right]\\
&-\frac{\rho}{2}\left\|d_{z}^{k}\right\|^{2}-\frac{\rho\beta_{x}}{2}\left\|\tilde{d}_{x}^{k}\right\|^{2}-\frac{\tau}{\rho}\left\|d_{y}^{k}\right\|^{2} - \delta \rho\|\hat{d}_x^k\|^2 \\
=& \frac{8\left(1+\tau\right)c^2_x\rho}{ \sigma^{\text{min}}_{A^{\top}A}}\left(\left\|\tilde{d}_{x}^{k-1}\right\|^{2}-\left\|\tilde{{d}}_{x}^{k}\right\|^{2}\right)\\
&+\frac{8\left(1+\tau\right)L^2}{\rho \sigma^{\text{min}}_{A^{\top}A}}\left(\left\|\hat{d}_{x}^{k-1}\right\|^{2}-\left\|\hat{{d}}_{x}^{k}\right\|^{2}\right)\\
&+\frac{8\left(1+\tau\right)\rho c^2_x}{ \sigma^{\text{min}}_{A^{\top}A}}\left(\|d_z^{k-1}\|^2-\|d_z^k\|^2\right)\\
&-\frac{\rho}{2}\left\|d_{z}^{k}\right\|^{2}-\frac{\rho\beta_{x}}{2}\left\|\tilde{d}_{x}^{k}\right\|^{2}-\frac{\tau}{\rho}\left\|d_{y}^{k}\right\|^{2} - \delta \rho\|\hat{d}_x^k\|^2\\
&+\frac{\left(1+\tau\right)\rho\left[2\left(\frac{L}{\rho}+c_x\right)^2+8c_x^2\right]}{ \sigma^{\text{min}}_{A^{\top}A}}\|\tilde{d}_x^k\|^2\\
&+\frac{8\left(1+\tau\right)L^2}{\rho\sigma^{\text{min}}_{A^{\top}A}}\|\hat{d}^k_x\|^2+\frac{16\left(1+\tau\right)\rho c_x^2}{\sigma^{\text{min}}_{A^{\top}A}}\|d_z^k\|^2\\
\end{aligned}
\end{equation}

\begin{equation*}
\begin{aligned}
=& \frac{8\left(1+\tau\right)c^2_x\rho}{ \sigma^{\text{min}}_{A^{\top}A}}\left(\left\|\tilde{d}_{x}^{k-1}\right\|^{2}-\left\|\tilde{{d}}_{x}^{k}\right\|^{2}\right)\\
&+\frac{8\left(1+\tau\right)L^2}{\rho \sigma^{\text{min}}_{A^{\top}A}}\left(\left\|\hat{d}_{x}^{k-1}\right\|^{2}-\left\|\hat{{d}}_{x}^{k}\right\|^{2}\right)\\
&+\frac{8\left(1+\tau\right)\rho c^2_x}{\sigma^{\text{min}}_{A^{\top}A}}\left(\|d_z^{k-1}\|^2-\|d_z^k\|^2\right)\\
&-\frac{\tau\rho\beta_{x}}{2}\left\|\tilde{d}_{x}^{k}\right\|^{2}- \tau\rho\|\hat{d}_x^k\|^2-\frac{\tau\rho}{2}\left\|d_{z}^{k}\right\|^{2}-\frac{\tau}{\rho}\left\|d_{y}^{k}\right\|^{2} \\
&-\left(1+\tau\right)\rho\left(\|\tilde{d}^k_x\|^2_{\tilde{\mathcal{D}}_x}+\|\hat{d}_x^k\|^2_{\hat{\mathcal{D}}_x}+\|d_z^k\|^2_{\mathcal{D}_z}\right)\\
&\leq \frac{8\left(1+\tau\right)c^2_x\rho}{ \sigma^{\text{min}}_{A^{\top}A}}\left(\left\|\tilde{d}_{x}^{k-1}\right\|^{2}-\left\|\tilde{{d}}_{x}^{k}\right\|^{2}\right)\\
&+\frac{8\left(1+\tau\right)L^2}{\rho \sigma^{\text{min}}_{A^{\top}A}}\left(\left\|\hat{d}_{x}^{k-1}\right\|^{2}-\left\|\hat{{d}}_{x}^{k}\right\|^{2}\right)\\
&+8\left(1+\tau\right)\rho \frac{c_x^2}{\sigma^{\min}_{A^{\top}A}} \left(\|d_z^{k-1}\|^2-\|d_z^k\|^2\right)\\
&-\frac{\tau\rho\beta_{x}}{2}\left\|\tilde{d}_{x}^{k}\right\|^{2}- \tau\rho\|\hat{d}_x^k\|^2-\frac{\tau\rho}{2}\left\|d_{z}^{k}\right\|^{2}-\frac{\tau}{\rho}\left\|d_{y}^{k}\right\|^{2}, \\
\end{aligned}
\end{equation*}
where $0<\tau<\delta<1$, and $\tilde{\mathcal{D}}_x\succeq 0$, $\hat{\mathcal{D}}_x\succeq 0$, and $\mathcal{D}_z\succeq 0$ are defined in \eqref{eq:D_x_tild}, \eqref{eq:D_x_hat}, and \eqref{eq:D_z}, respectively.
Thus, \eqref{eq:E_ineq} follows from \eqref{eq:theo_1_proof_3} and the definition of $E^{k+1}$ in \eqref{eq:E_tild}.
Similarly, by \eqref{cond:alpha_subp_cond} and $\hat{d}^k_x=x^{k+1}-\tilde{x}^{k+1}$, we have
\begin{equation}
    \label{eq:L_w_tild_ineq}
    \begin{array}{rl}
        &\mathcal{L}_{\rho}(\tilde{x}^{k+1}, z^{k+1}, y^{k+1}) - \mathcal{L}_{\rho}(\tilde{x}^{k}, z^{k}, y^{k}) \\
        \leq& \mathcal{L}_{\rho}(\tilde{x}^{k+1}, z^{k+1}, y^{k+1}) - \mathcal{L}_{\rho}(x^{k}, z^{k}, y^{k})-\delta\rho \|x^k-\tilde{x}^k\|^2 \\
        =&\mathcal{L}_{\rho}(\tilde{x}^{k+1}, z^{k+1}, y^{k+1}) - \mathcal{L}_{\rho}(x^{k}, z^{k}, y^{k})-\delta\rho \|\tilde{d}^{k-1}_{x}\|^2
    \end{array}
\end{equation}
Using the definition of $\tilde{E}^{k+1}$ in \eqref{eq:E_tild}, the same argument yields \eqref{eq:E_tild_ineq}.

\subsection{Proof of Theorem~\ref{theo:convergence}}
\label{app:proof_theorem}
Assume that the energy sequence $\{\tilde{E}^k\}$ is bounded from below.
From the energy descent inequality \eqref{eq:E_tild_ineq}, we obtain
\begin{equation}
    \label{eq:E_lower_bound_ineq}
    \begin{array}{cl}
         &c\sum_{k=1}^{K}\left\{\left\|\tilde{d}_{x}^{k}\right\|^{2}+\|\hat{d}_x^k\|^2+\left\|d_{z}^{k}\right\|^{2}+\left\|d_{y}^{k}\right\|^{2}\right\} \\
         & \leq \tilde{E}^1-\tilde{E}^{K+1} \leq \tilde{E}^1-\bar{P}
    \end{array}
\end{equation}
where $c = \operatorname{min}\left\{ \frac{\tau\rho\beta_{x}}{2}, \frac{\tau\rho}{2}, \tau\rho, \frac{\tau}{\rho}\right\}>0$ and $\bar{P}$ denotes the lower bound of $\tilde{E}^k$. 
This implies that the corresponding increments are square summable, and hence
\begin{equation}
    \label{eq:d_limit}
    \begin{split}
        \lim _{k \rightarrow \infty}\left\|\tilde{d}_{x}^{k}\right\|=0, \quad \lim _{k \rightarrow \infty}\left\|\hat{d}_{x}^{k}\right\|=0, \\
        \lim _{k \rightarrow \infty}\left\|d_{z}^{k}\right\|=0 \quad \text { and } \quad \lim _{k \rightarrow \infty}\left\|d_{y}^{k}\right\|=0.
    \end{split}
\end{equation}
Using $d_y^k=\rho \tilde{r}^{k+1}$ and the definition of $R^k$ in \eqref{eq:R_k+1}, we further obtain
\begin{equation}
    \label{eq:limit_r_R}
    \begin{split}
        &\lim _{k \rightarrow \infty}\left\|\tilde{r}^{k}\right\|=0 \quad \text { and } \\
        &\lim _{k \rightarrow \infty} R^{k}=\lim _{k \rightarrow \infty}\left(\left\|\tilde{d}_{x}^{k-1}\right\|+\left\|d_{z}^{k-1}\right\|+\left\|\tilde{r}^{k}\right\|\right)=0.
    \end{split}
\end{equation}
The residual decomposition $r^k=\tilde{r}^k+A\hat{d}_{x}^{k-1}$, together with the primal update bound $\|d_x^k\|\leq \|\hat{d}_x^k\|+\|\tilde{d}_x^k\|$, gives
\begin{equation}
    \label{eq:limit_r_d}
    \lim _{k \rightarrow \infty}\left\|r^{k}\right\|=0 \quad \text { and } \quad \lim _{k \rightarrow \infty}\left\|d_{x}^{k}\right\|=0
\end{equation}
where $r^k=Ax^k-z^k$ and $d^k_x=x^{k+1}-x^k$.
Since $\{\tilde{E}^k\}$ is monotonically decreasing and bounded from below, there exists $F^{*}$ such that $\lim _{k \rightarrow \infty} \tilde{E}^{k}=F^{*}$.
Then, by the definition of $E^k$ together with \eqref{eq:limit_r_R} and \eqref{eq:limit_r_d}, we obtain \eqref{eq:lagran_limit}.

We next establish stationarity by examining the Lagrangian subdifferentials:
\begin{equation}
    \label{eq:partial_L_rho}
    \begin{aligned}
\partial_{x} \mathcal{L}_{\rho}\left(x^k, z^k, y^k\right) & =\partial_{x} \mathcal{L}\left(x^k, z^k, y^k\right)+\rho A^{\top} r^{k}\\
&=\nabla f\left(x^{k}\right)+A^{\top} y^{k}+\rho A^{\top} r^{k} \\
& =\nabla_{x} \mathcal{L}_{\rho}\left(\tilde{x}^{k}, z^{k}, y^{k-1}\right)+A^{\top}d_y^{k-1}\\
&+\rho A^{\top} A \hat{d}_{x}^{k-1}+\left[\nabla f\left(x^k\right)-\nabla f\left(\tilde{x}^{k}\right)\right], \\
\partial_{z} \mathcal{L}_{\rho}\left(x^k, z^k, y^k\right) & =\partial_{z} \mathcal{L}\left(x^k, z^k, y^k\right)-\rho r^{k}\\
&=\partial_{z} g\left(z^{k}\right)-y^{k}-\rho r^{k} \\
& =\partial_{z} \mathcal{L}_{\rho}\left({x}^{k-1}, z^{k}, y^{k-1}\right)-d_y^{k-1}-\rho Ad^{k-1}_x, \\
\partial_{y} \mathcal{L}_{\rho}\left(x^k, z^k, y^k\right) & =\partial_{y} \mathcal{L}\left(x^k, z^k, y^k\right)=r^{k} .
\end{aligned}
\end{equation}
Since subproblem \eqref{prob:inexact_z_sub_prob} is solved exactly,
\begin{equation}
\label{eq:part_z_subprob}
0 \in \partial_{z} \mathcal{L}_{\rho}\left({x}^{k-1}, z^{k}, y^{k-1}\right).
\end{equation}
Together with condition \eqref{cond:x_subp_cond_2}, \eqref{eq:d_limit}, and \eqref{eq:limit_r_d}, this yields \eqref{eq:paritial_L_limit}.
For any limiting point $(x^{*}, z^{*}, y^{*})$ of $\left\{x^k, z^k, y^k\right\}$, \eqref{eq:paritial_L_limit} and the definition of the limiting subdifferential $\partial \mathcal{L}(x^{*}, z^{*}, y^{*})$ imply that \eqref{cond:optimal_conds} holds. Therefore, $(x^{*}, z^{*}, y^{*})$ is a stationary point of \eqref{prob:convex_qp_3}.
Since the QPs considered in this work are convex, these stationary points further guarantee global optimality.

Finally, \eqref{eq:E_lower_bound_ineq} implies the following bound on the minimum residual:
\begin{equation}
    \begin{aligned}
        \min _{k \in\{1, \ldots, K\}}\left\{\left\|\hat{d}_{x}^{k-1}\right\|^{2}+\left\|\tilde{d}_{x}^{k}\right\|^{2}+\left\|d_{z}^{k}\right\|^{2}+\left\|\tilde{r}^{k+1}\right\|^{2}\right\}\\
        =\mathcal{O}(1 / K),
    \end{aligned}
\end{equation}
Together with \eqref{cond:x_subp_cond_2} and \eqref{eq:part_z_subprob}, this implies
\begin{equation}
    \min _{k \in\{1, \ldots, K\}}\left\{\operatorname{dist}\left(0, \partial \mathcal{L}\left(x^k, z^k, y^k\right)\right)\right\}=\mathcal{O}(1 / \sqrt{K}).
\end{equation}
This completes the proof that Algorithm~\ref{alg:osqp_lstm} converges sublinearly to stationary points under the stated assumptions, with the rate governed by the problem geometry.


\end{document}